# A DIFFUSION MODEL OF SCHEDULING CONTROL IN QUEUEING SYSTEMS WITH MANY SERVERS[1]

### By Rami Atar

*Technion—Israel Institute of Technology*

This paper studies a diffusion model that arises as the limit of a queueing system scheduling problem in the asymptotic heavy traffic regime of Halfin and Whitt. The queueing system consists of several customer classes and many servers working in parallel, grouped in several stations. Servers in different stations offer service to customers of each class at possibly different rates. The control corresponds to selecting what customer class each server serves at each time. The diffusion control problem does not seem to have explicit solutions and therefore a characterization of optimal solutions via the Hamilton–Jacobi–Bellman equation is addressed. Our main result is the existence and uniqueness of solutions of the equation. Since the model is set on an unbounded domain and the cost per unit time is unbounded, the analysis requires estimates on the state process that are subexponential in the time variable. In establishing these estimates, a key role is played by an integral formula that relates queue length and idle time processes, which may be of independent interest.

**1. Introduction.** We consider optimal scheduling control for a class of queueing systems that operate in heavy traffic, in the sense that the load on the system is nearly equal to its capacity. As often occurs, exact analysis of the control problem is unavailable and an asymptotic approach is taken, where a parametrization of the model is introduced and a diffusion control problem is obtained in the limit. The parametrization that has been more common in research papers on related problems (referred to here as conventional heavy traffic) is one where arrival and service rates are both scaled up

---
Received August 2003; revised April 2004.

[1]Supported in part by the Israel Science Foundation Grant 126/02, the US–Israel Binational Science Foundation Grant 1999179 and the fund for promotion of research at the Technion.

*AMS 2000 subject classifications.* 60K25, 68M20, 90B22, 90B36, 49L20.

*Key words and phrases.* Multiclass queueing systems, scheduling and routing control, heavy traffic, Halfin–Whitt regime, buffer-station tree, control of diffusions, Hamilton–Jacobi–Bellman equation, unbounded domain.

in a way that the system operates near full capacity. Recently, several papers have studied a different parametrization, proposed by Halfin and Whitt [10], where increase of arrivals is balanced by scaling up the *number* of (identical) servers in each service station, while keeping the service time distribution of the individual servers fixed. In the limit as the parameter grows without bound, conventional heavy traffic typically gives rise to diffusion processes in the orthant with reflection on the boundary, whereas in the Halfin–Whitt (HW) regime the diffusion takes values in the Euclidean space. This paper focuses on the diffusion model that corresponds to the queueing system introduced below, operating in the HW regime.

The queueing system has a fixed number of customer classes and many exponential servers grouped in a fixed number of stations. Only some stations can offer service to each class, and the service rates depend on the class and on the station. Also, customers not being served may abandon the system (see Figure 1). Scheduling (and routing) of jobs in the queueing system is regarded as control. As cost one considers an expected discounted cumulative function of performance criteria such as queue lengths, number of idle servers or number of customers of each class present at each station. The system is parametrized so that the arrival rates and the number of servers at each station are nearly proportional to a large parameter $n$, while service and abandonment rates are nearly constant. For motivation on the model and on this asymptotic regime, see [4] and [12].

In the scaling limit and under appropriate assumptions, one obtains a diffusion model whose ingredients are denoted by $X, Y, Z, \Psi$ and $\widetilde{W}$. Let $\mathcal{I}$ and $\mathcal{J}$ be index sets for customer classes and service stations, respectively. Assume that a sequence of systems is given where, in the $n$th system, the number of servers at each station is proportional to $n$. For each system $n$, denote by $X_i^n, i \in \mathcal{I}$, the number of customers of class $i$ present in the system. Then $X$ stands for the (formal) weak limit of the processes $n^{-1/2}(X^n - nx^*)$, where $x^*$ are constants that come from a corresponding "static fluid model" (see [2] for full details on the fluid model as well as on the derivation of the diffusion model from the queueing model). Note that because of the centering, $X_i$ assume both positive and negative values. In an analogous fashion the processes $Y_i$, $i \in \mathcal{I}$, and $Z_j$, $j \in \mathcal{J}$, correspond to queue length of class $i$ and, respectively, the number of idle servers at station $j$, and $\Psi_{ij}$ corresponds to the number of class-$i$ customers in service at station $j$. Finally, $\widetilde{W}_i$ are Brownian motions that represent the effect of fluctuations in arrival and service times. With $\theta_i$ and $\mu_{ij}$ standing for the abandonment rate of class-$i$ customers and, respectively, service rate of class-$i$ customers at station $j$, the diffusion model reads

$$(1) \qquad X_i(t) = x_i + \widetilde{W}_i(t) - \sum_{j \in \mathcal{J}} \mu_{ij} \int_0^t \Psi_{ij}(s)\,ds - \theta_i \int_0^t Y_i(s)\,ds, \qquad i \in \mathcal{I},$$

$$(2) \qquad \sum_{j \in \mathcal{J}} \Psi_{ij} = X_i - Y_i, \qquad i \in \mathcal{I},$$

$$(3) \qquad \sum_{i \in \mathcal{I}} \Psi_{ij} = -Z_j, \qquad j \in \mathcal{J},$$

$$(4) \qquad Y_i \geq 0, \qquad Z_j \geq 0, \qquad i \in \mathcal{I}, j \in \mathcal{J}.$$

$$(5) \qquad \min\left[\sum_{i \in \mathcal{I}} Y_i, \sum_{j \in \mathcal{J}} Z_j\right] = 0.$$

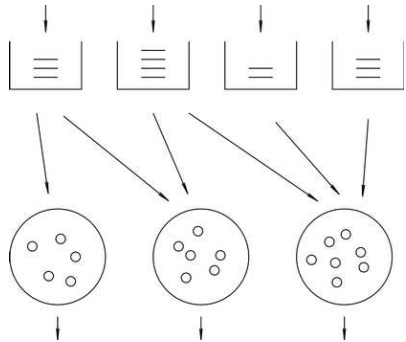

Fig. 1. *System with four customer classes (buffers) and three server types (service stations).*

Denote by $\mathcal{T}$ the graph that has a node per each class and each station, and an edge that joins the class and the station if and only if customers of the class can be served at the station. It is assumed in this paper that $\mathcal{T}$ is a tree (see Section 2 for discussion on various aspects of this assumption). It is then possible to encode equations (1)–(5) in a single equation of the form $X(t) = x + rW(t) + \int_0^t b(X(s), U(s))\,ds$, where $U$ is a control process that takes values in a compact space and $W$ is a standard Brownian motion. Also, rescaling the cost associated with the queueing system appropriately results in a cost of the form $E_x \int_0^\infty e^{-\gamma t} L(X_t, U_t)\,dt$ for the diffusion model.

The main result of this paper is characterization of the diffusion control problem's value as the unique solution to the associated Hamilton–Jacobi–Bellman (HJB) equation. Such problems for diffusion in an unbounded domain are well understood when either the "running cost" $L$ is bounded [5, 7] or the drift $b$ is bounded (as follows from the results of [13], e.g.). The difficulty in the current model stems from the fact that the domain as well as the functions $L$ and $b$ are not necessarily bounded. The question then requires deeper understanding of the model and, in particular, estimates on moments of the process that are subexponential in the time variable. Our results apply when either one of the following conditions holds (further conditions on abandonment rates are assumed in each case; see Section 2): 1. Service rates are either class- or station-dependent; 2. $\mathcal{T}$ satisfies diam($\mathcal{T}$) $\leq 3$; 3. Running cost $L(X, U)$ is comparable to $\|X\|^m$, some $m \geq 1$; 4. $L$ is bounded.

The large time estimates needed for the first three cases of the main result use three different strategies. One of our basic tools is an integral equation (51) that expresses a relationship directly between $Y$, $Z$ and $\widetilde{W}$, not involving $X$ and $\Psi$. The special form this equation takes when $\mu_{ij}$ depend only on $i$ or only on $j$ makes it possible to get the required estimates on $X$ in case 1 of the main result. Equation (51) may be of independent interest since it can be seen to express a relationship between the data $\widetilde{W}$, the control process $U$ and the *one-dimensional* process $\sum_i X_i$ alone (cf. Section 4).

The integral equation is used also in case 2 of the main result, along with a certain property of the system (1)–(5). This property, which we call the *nonidling property*, should be understood as one of the system of equations rather than the stochastic processes, because it does not regard $\widetilde{W}$ as a Brownian motion, but as a generic function:

> *If the system starts with $x_i > 0$ and $t \mapsto \widetilde{W}_i(t)$ are strictly increasing, then $Z_j(t) = 0$, $j \in \mathcal{J}$, $t \geq 0$.*

For trees of diameter 3 or less we can show that the property holds and that it implies moment estimates that are polynomial in time. Heuristically this property expresses that servers do not idle when there is large enough amount of work. Viewed this way it appears to be a basic question on the model and one would like to understand in what generality it holds (e.g., for more general trees), irrespective of the goals of the current work.

In case 3 of the main result we assume that the function $L$ is comparable to a power of the norm of the state process $X$. This assumption is not the most natural in the context of queueing, since, for example, cost functions that depend on queue lengths cannot be treated. However, it can be useful if one is interested in stabilizing the dynamical system about a nominal model, since the cost penalizes deviations from the static fluid model. On the technical side, penalizing deviations from a nominal model simplifies the problem in that moment estimates are required under one particular control rather than under all controls, and as a result we can treat the model at full generality, as far as the tree structure and the service rates are concerned. Finally, in case 4, where $L$ is assumed to be bounded, estimates on $X$ are not needed.

Apart from their own contribution the results of this paper are a first step toward identifying scheduling policies for the queueing system that are, in an appropriate sense, asymptotically optimal [2]. As in [4], such policies can be derived from the solution to the diffusion control problem. In addition, such asymptotic analysis justifies the relationship between the queueing system and the diffusion model studied here.

Recent results on the HW regime include the following. Puhalskii and Reiman [17] extended the work of Halfin and Whitt to multiple customer classes, priorities and phase-type service distribution. Mandelbaum, Massey and Reiman [15] established functional law of large numbers and central limit theorems for a wide class of Markovian systems in the HW regime. Armony and Maglaras [1] modeled and analyzed rational customers in equilibrium, and Garnett, Mandelbaum and Reiman [8] studied models with abandonments from the queue. Papers where a control theoretic approach was taken to study queueing systems in this regime are few. The diffusion control problem associated with scheduling jobs in a system with multiple customer classes and a single service station was analyzed by Harrison and Zeevi [12]. In a similar setting, Atar, Mandelbaum and Reiman [4] established asymptotic optimality of scheduling policies for the queueing system derived from the diffusion model's HJB equation. A special case where explicit, pathwise solutions to the diffusion control problem are available appeared in [3]. Finally, with regard to the conventional heavy traffic analysis of systems analogous to the one studied in the current paper, we mention Harrison and López [11], where the corresponding Brownian control problem is identified and solved, Williams [18], where a dynamic threshold scheduling policy is proposed for the queueing system, and Mandelbaum and Stolyar [16], where asymptotic optimality of a simple scheduling policy is proved for convex delay costs.

The organization of the paper is as follows. Section 2 describes the queueing system and the diffusion model and states the main result. Section 3 shows how estimates on $\|X\|$ reduce to estimates on $|\sum_i X_i|$. Section 4 develops the integral equation (51) and establishes moment estimates in case 1 of the main result. Section 5 studies the nonidling property for trees of diameter 3 and establishes case 2 of the main result. Section 6 treats case 3 and summarizes the estimates in all cases. Based on the estimates of Sections 3–6, the proof of the main result is similar to an analogous treatment in [4], but for completeness we provide it in the Appendix.

NOTATION. For $f : [0, \infty) \to \mathbb{R}$, let $\Im f = \int_0^{\cdot} f$. For a vector $x$, let $\|x\| = \sum |x_i|$. For two column vectors $v$ and $u$, $v \cdot u$ denotes their scalar product. The symbols $e_i$ denote the coordinate unit vectors and $e = (1, \ldots, 1)'$. The dimension of $e$ may change from one expression to another, and, for example, $e \cdot a = e \cdot b$ makes sense even if $a$ and $b$ are of different dimension. The symbol $C^{m,\varepsilon}$ (resp. $C^m$) denotes the class of functions on $\mathbb{R}^I$ for which all derivatives up to order $m$ are Hölder continuous uniformly on compacts (resp. continuous); $C_{\mathrm{pol}}$ denotes the class of continuous functions $f$ on $\mathbb{R}^I$, for which there are $c$, $r$ such that $|f(x)| \le c(1 + \|x\|^r)$, $x \in \mathbb{R}^I$. Let $C_{\mathrm{pol}}^{m,\varepsilon} = C_{\mathrm{pol}} \cap C^{m,\varepsilon}$ and $C_{\mathrm{pol}}^m = C_{\mathrm{pol}} \cap C^m$, and let $C_{\mathrm{pol},+}^m$ be the class of nonnegative functions in $C_{\mathrm{pol}}^m$. The symbol $C_b$ denotes the class of bounded

continuous functions and $C_b^m = C^m \cap C_b$. Let $\mathbb{R}_+ = [0, \infty)$. If $X$ is a process or a function on $\mathbb{R}_+$, $\|X\|_t^* = \sup_{0 \le s \le t} \|X(s)\|$, and if $X$ takes real values, $|X|_t^* = \sup_{0 \le s \le t} |X(s)|$. The symbols $X(t)$ and $X_t$ are used interchangeably. The symbols $c_1, c_1', c_2, c_2', \ldots$ denote deterministic positive constants that may have different values in the proof of different results.

**2. The queueing system and the diffusion.** We start with an example that demonstrates how the HW scaling is performed on a simple queueing system. The queueing system, parametrized by $n \in \mathbb{N}$, has a single customer class with renewal arrivals at rate $\lambda^n$ and $n$ servers, each having exponential service time distribution of rate $\mu^n$. While a customer is not in service, it abandons the system at rate $\theta^n$. Let $X^n(t)$, $Y^n(t)$ and $Z^n(t)$ denote the total number of customers, the number of customers not being served and the number of servers that are idle at time $t$, respectively. Clearly

$$(6) \qquad\qquad X^n + Z^n = n + Y^n.$$

Assume that the system operates under work conservation, in the sense that $Y^n \wedge Z^n = 0$. Then $Y^n = (X^n - n)^+$ and $Z^n = (X^n - n)^-$. The parameters scale so that $n^{1/2}(n^{-1}\lambda^n - \lambda) \to \hat{\lambda}$, $n^{1/2}(\mu^n - \mu) \to \hat{\mu}$ and $\theta^n \to \theta$. The system is assumed to be critically loaded, in the sense that

$$(7) \qquad\qquad \lambda = \mu.$$

Denote $\hat{X}^n(t) = n^{-1/2}(X^n(t) - n)$, $\hat{Y}^n(t) = n^{-1/2}Y^n(t)$ and $\hat{Z}^n(t) = n^{-1/2}Z^n(t)$. Assuming that the interarrival times have finite second moment and that $\hat{X}^n(0) \to x$, $X^n$ converges weakly to a diffusion $X$ that solves

$$(8) \qquad X(t) = x + rW(t) + \int_0^t (\hat{\lambda} - \hat{\mu} + \mu X(s)^- - \theta X(s)^+)\, ds,$$

where $W$ is a standard Brownian motion and $r$ is a constant that depends on the first two moments of the interarrival time and on $\lambda$ (in the absence of abandonment, the result is due to [10]; see [8] for a treatment of abandonment). Denoting by $Y$ and $Z$ the weak limit of $\hat{Y}^n$ and $\hat{Z}^n$, respectively, note that (8) can be rewritten as

$$(9) \qquad X(t) = x + rW(t) + \int_0^t (\hat{\lambda} - \hat{\mu} + \mu Z(s) - \theta Y(s))\, ds,$$

$$(10) \qquad\qquad X = Y - Z, \qquad Y \wedge Z = 0.$$

We repeat that some nonexponential service distributions were treated in [17]; however, the diffusion limit turns out to be more complicated and that approach is not taken here.

The diffusion model studied in this paper corresponds to a queueing system with several classes and stations. The queueing system and the corresponding diffusion model are introduced below. Since in the current paper we focus on the diffusion model, we present it here without attempting to justify its relationship to the queueing system, and we demonstrate only how it is analogous to the simple model (8)–(10). Full details on this relationship are deferred to [2]. The queueing system has $I$ customer classes and $J$ service stations (see Figure 1). At each station there are many independent servers of the same type. Each customer requires service only once and can be served indifferently by any server at the same station, but possibly at different rates at different stations. Only some stations can offer service to each class. Label the classes (and corresponding buffers) as $1, \ldots, I$ and the types (and corresponding stations) as $I + 1, \ldots, I + J$, and set

$$\mathcal{I} = \{1, \ldots, I\}, \qquad \mathcal{J} = \{I + 1, \ldots, I + J\}.$$

The structure of the system can be encoded in a graph. A pair $(i, j) \in \mathcal{I} \times \mathcal{J}$ is called an *activity* if customers of class $i$ can be served at station $j$. It is assumed that whether $(i, j)$ is an activity does not depend on $n$. Let $\mathcal{T}$ denote the graph with vertex set $\{1, 2, \ldots, I + J\} = \mathcal{I} \cup \mathcal{J}$: A node is associated with each buffer and each station. Edges of $\mathcal{T}$ are between elements $i \in \mathcal{I}$ and $j \in \mathcal{J}$ such that $(i, j)$ is an activity. Write $i \sim j$ if $(i, j)$ is an activity. Denote the edge set for the graph by

$$\mathcal{E} = \{(i, j) \in \mathcal{I} \times \mathcal{J} : i \sim j\}.$$

For $j \in \mathcal{J}$, let $N_j^n$ be the number of servers at station $j$. Let $X_i^n(t)$, $Y_i^n(t)$ and $Z_j^n(t)$ denote the total number of class-$i$ customers in the system, the number of class-$i$ customers in the queue, and the number of idle servers in station $j$ at time $t$, respectively. Finally, let $\Psi_{ij}^n(t)$ denote the number of class-$i$ customers in service at station $j$ at time $t$ (note that $\Psi_{ij}^n = 0$, $i \nsim j$). In vector–matrix notation, set $X^n = (X_i^n)_{i \in \mathcal{I}}$, $Y^n = (Y_i^n)_{i \in \mathcal{I}}$, $Z^n = (Z_j^n)_{j \in \mathcal{J}}$ and $\Psi^n = (\Psi_{ij}^n)_{i \in \mathcal{I}, i \in \mathcal{J}}$. Straightforward relationships are expressed by the equations

$$(11) \qquad Y_i^n + \sum_{j \in \mathcal{J}} \Psi_{ij}^n = X_i^n, \qquad i \in \mathcal{I},$$

$$(12) \qquad Z_j^n + \sum_{i \in \mathcal{I}} \Psi_{ij}^n = N_j^n, \qquad j \in \mathcal{J},$$

$$(13) \qquad Y_i^n(t), Z_j^n(t) \geq 0, \qquad i \in \mathcal{I}, j \in \mathcal{J}, t \geq 0.$$

Arrivals of class-$i$ customers occur at rate $\lambda_i^n$, abandonment from queue $i$ is at rate $\theta_i^n$ and, for $i \sim j$, service of a class-$i$ customer at station $j$ is at rate $\mu_{ij}^n$. Assume

$$(14) \qquad n^{1/2}(n^{-1}\lambda_i^n - \lambda_i) \to \hat{\lambda}_i, \qquad n^{1/2}(\mu_{ij}^n - \mu_{ij}) \to \hat{\mu}_{ij}, \qquad \theta_i^n \to \theta_i,$$

$$(15) \qquad n^{1/2}(n^{-1}N_j^n - \nu_j) \to 0.$$

Above, $\lambda_i, \nu_j \in (0, \infty)$, $\theta_i \in [0, \infty)$, $\hat{\lambda} \in \mathbb{R}$ and for $i \sim j$, $\mu_{ij} \in (0, \infty)$ and $\hat{\mu}_{ij} \in \mathbb{R}$. It is convenient to set throughout $\mu_{ij} = \hat{\mu}_{ij} = 0$ for $i \nsim j$. In analogy with (7) in the simple example, we invoke a condition that expresses that the system is critically loaded. The condition involves the "first order" parameters $\lambda_i, \nu_j$ and $\mu_{ij}$, and certain constants $x_i^*$ and $\psi_{ij}^*$ that represent a static fluid model (see details in [2]). The processes are rescaled as

$$(16) \qquad \widehat{X}_i^n(t) = n^{-1/2}(X_i^n(t) - nx_i^*),$$

$$(17) \qquad \widehat{Y}_i^n(t) = n^{-1/2}Y_i^n(t), \qquad \widehat{Z}_j^n(t) = n^{-1/2}Z_j^n(t),$$

$$\widehat{\Psi}_{ij}^n(t) = n^{-1/2}(\Psi_{ij}^n(t) - n\psi_{ij}^*).$$

Assuming $\widehat{X}^n(0) \to x$ and considering $X$, $Y$, $Z$ and $\Psi$ as formal limits of $\widehat{X}^n$, $\widehat{Y}^n$, $\widehat{Z}^n$ and $\widehat{\Psi}^n$, we expect in analogy with (9),

$$(18) \quad X_i(t) = x_i + \widetilde{W}_i(t) - \sum_{j \in \mathcal{J}} \mu_{ij} \int_0^t \Psi_{ij}(s)\,ds - \theta_i \int_0^t Y_i(s)\,ds, \qquad i \in \mathcal{I},$$

holds, where $\widetilde{W}_i(t) = r_i W_i(t) + \ell_i t$, $W$ is a standard Brownian motion, $r_i \in (0, \infty)$ and $\ell_i \in \mathbb{R}$ are constants, $\Psi_{ij} = 0$ for $i \nsim j$, and, in view of (11)–(13),

$$(19) \qquad\qquad \sum_{j \in \mathcal{J}} \Psi_{ij} = X_i - Y_i, \qquad i \in \mathcal{I},$$

$$(20) \qquad\qquad \sum_{i \in \mathcal{I}} \Psi_{ij} = -Z_j, \qquad j \in \mathcal{J},$$

$$(21) \qquad\qquad Y_i \geq 0, \qquad Z_j \geq 0, \qquad i \in \mathcal{I}, j \in \mathcal{J}.$$

For reasons explained in [2], the work conservation condition is replaced by the condition

$$(22) \qquad\qquad e \cdot Y \wedge e \cdot Z = 0.$$

The diffusion model is now described by (18)–(22). To view the model in a control theoretic setting, regard $\Psi$ as a control process and $X$ as a controlled diffusion. Then (18) describes the dynamics of $X$, and (19)–(22) serve to define the "auxiliary" processes $Y$ and $Z$ and to set constraints on $\Psi$. Note that the constraints on $\Psi$ involve the process $X$.

While relationships (18)–(22) are, in a sense, analogous to (9)–(10) (although obviously there is no control process in the simple model), we would like also to have a relationship analogous to (8). More precisely, we seek to describe the model in the convenient form

$$(23) \qquad\qquad X(t) = x + rW(t) + \int_0^t b(X(s), U(s))\,ds$$

with $r = \mathrm{diag}(r_i)_{i \in \mathcal{I}}$, appropriate function $b$ and control process $U$, and, in particular, without having to impose constraints on $U$ that involve $X$. The assumption below is useful in this development; however, as discussed at the end of this section, the prime reason for imposing it is different.

ASSUMPTION 1 (Treelike). *The graph $\mathcal{T}$ is a tree.*

Proposition A.2 in the Appendix shows that under Assumption 1, (19) and (20) are equivalent to

$$\Psi = G(X - Y, -Z), \tag{24}$$

$G$ being a linear map from $\{(\alpha, \beta) \in \mathbb{R}^{I+J} : \sum \alpha_i = \sum \beta_j\}$ to $\mathbb{R}^{IJ}$. We proceed to derive (23) under the treelike assumption. Note first that by (19) and (20), $e \cdot X = e \cdot Y - e \cdot Z$, and thus by (21) and (22), $e \cdot Y = (e \cdot X)^+$ and $e \cdot Z = (e \cdot X)^-$. Hence $Y$ can be represented as

$$Y_i(t) = (e \cdot X(t))^+ u_i(t), \tag{25}$$

where $u_i(t) \geq 0$ and $e \cdot u(t) = 1$. Similarly,

$$Z_j(t) = (e \cdot X(t))^- v_j(t), \tag{26}$$

where $v_j(t) \geq 0$ and $e \cdot v(t) = 1$. Consider $U := (u, v)$ as a control process that takes values in

$$\mathbb{U} := \{(u, v) \in \mathbb{R}^{I+J} : u_i, v_j \geq 0, \ e \cdot u = e \cdot v = 1\}.$$

By (24),

$$\Psi = \widehat{G}(X, U) := G(X - (e \cdot X)^+ u, -(e \cdot X)^- v).$$

Let

$$b_i(X, U) = -\sum_{j \in \mathcal{J}} \mu_{ij} \widehat{G}(X, U)_{ij} - \theta_i (e \cdot X)^+ u_i + \ell_i \tag{27}$$

and write $b = (b_i)_{i \in \mathcal{I}}$. We see that (18) can be written as (23).

DEFINITION 1 (Admissible systems and controlled processes).

(i) We call

$$\pi = (\Omega, F, (F_t), P, U, W)$$

an *admissible system* and say that $U$ is a *control* associated with $\pi$ if $(\Omega, F, (F_t), P)$ is a complete filtered probability space, $U$ is an $(F_t)$ progressively measurable $\mathbb{U}$-valued process and $W$ is a standard $I$-dimensional $(F_t)$ Brownian motion.

(ii) We say that $X$ is a *controlled process* associated with initial data $x \in \mathbb{R}^I$ and an admissible system $\pi = (\Omega, F, (F_t), P, U, W)$ if $X$ is a continuous, $(F_t)$-adapted process such that $P$-a.s., $\int_0^t |b(X(s), U(s))|\, ds < \infty$ and

$$(28) \qquad X(t) = x + rW(t) + \int_0^t b(X(s), U(s))\, ds, \qquad 0 \le t < \infty.$$

As stated in Proposition A.1 in the Appendix, there is a unique controlled process $X$ associated with any $x$ and $\pi$. With an abuse of notation we sometimes denote the dependence on $x$ and $\pi$ by writing $P_x^\pi$ in place of $P$ and $E_x^\pi$ in place of $E$. Denote by $\Pi$ the class of all admissible systems.

Let a constant $\gamma > 0$ and a function $L$ be given, and consider the cost

$$C(x, \pi) = E_x^\pi \int_0^\infty e^{-\gamma t} L(X(t), U(t))\, dt, \qquad x \in \mathbb{R}^I, \pi \in \Pi.$$

Our assumption on $L$ is as follows.

Assumption 2. (i) We have $L(x, U) \ge 0$ and $(x, U) \in \mathbb{R}^I \times \mathbb{U}$.

(ii) The mapping $(x, U) \mapsto L(x, U)$ is continuous.

(iii) There is $\varrho \in (0, 1)$ such that for any compact $A \subset \mathbb{R}^I$,

$$|L(x, U) - L(y, U)| \le c\|x - y\|^\varrho$$

holds for $U \in \mathbb{U}$ and $x, y \in A$, where $c$ depends only on $A$.

(iv) There are constants $c_L > 0$ and $m_L \ge 1$ such that $L(x, U) \le c_L(1 + \|x\|^{m_L})$, $U \in \mathbb{U}$, $x \in \mathbb{R}^I$.

Define the value function as $V(x) = \inf_{\pi \in \Pi} C(x, \pi)$. The HJB equation for the problem is (cf. [7])

$$(29) \qquad\qquad \mathcal{L}f + H(x, Df) - \gamma f = 0,$$

where $\mathcal{L} = (1/2)\sum_{i \in \mathcal{I}} r_i^2 \partial^2 / \partial x_i^2$ and

$$H(x, p) = \inf_{U \in \mathbb{U}} [b(x, U) \cdot p + L(x, U)].$$

The equation is considered on $\mathbb{R}^I$ with the growth condition

$$(30) \qquad\qquad \exists c, m, \qquad |f(x)| \le c(1 + \|x\|^m), \qquad x \in \mathbb{R}^I.$$

Definition 2. Let $x \in \mathbb{R}^I$ be given. We say that a measurable function $h \colon \mathbb{R}^I \to \mathbb{U}$ is a *Markov control policy* if there is an admissible system $\pi$ and a controlled process $X$ corresponding to $x$ and $\pi$, such that $U_s = h(X_s)$, $s \ge 0$, $P$-a.s. We say that an admissible system $\pi$ is *optimal* for $x$ if $V(x) = C(x, \pi)$. We say that a Markov control policy is optimal for $x$ if at least one of the admissible systems corresponding to it is optimal.

Different parts of our main result below require different assumptions on the abandonment rates:

$$(31) \qquad \forall\, (i,j) \in \mathcal{E}, \qquad \theta_i \leq \mu_{ij},$$

$$(32) \qquad \exists\, (i,j) \in \mathcal{E}, \qquad \theta_i \leq \mu_{ij}.$$

These assumptions are rational for the following reason. If (31) does not hold [and certainly if (32) does not hold], there is a class $i$ where customers leave the system by abandonment more quickly than they do by getting served at a certain station $j \sim i$. Thus under many reasonable performance criteria (e.g., any increasing functional of the queue lengths) it is preferable to never use activity $(i,j)$. This stands in contrast to our work conservation assumptions.

**THEOREM 1.** *Let Assumptions* 1 *and* 2 *hold. In addition, let one of the following conditions hold.*

(i) *For* $(i,j) \in \mathcal{E}$, $\mu_{ij}$ *depends only on* $i$, *or for* $(i,j) \in \mathcal{E}$, $\mu_{ij}$ *depends only on* $j$. *In addition* $\theta_i = 0$, $i \in \mathcal{I}$.

(ii) *The tree* $\mathcal{T}$ *is of diameter* 3 *at most and* (31) *holds.*

(iii) *There are* $a_1, a_2 > 0$ *such that* $L(x,U) \geq a_1 \|x\|^{m_L}$ *for all* $\|x\| > a_2$ *and all* $U \in \mathbb{U}$ *(where* $m_L$ *is as in Assumption* 2*). In addition,* (32) *holds.*

(iv) *The function* $L$ *is bounded.*

*Then the value* $V$ *is in* $C^{2,\rho}_{\mathrm{pol}}$ *and it solves* (29) *and* (30). *In cases* (i) *and* (ii) *[resp.* (iii) *and* (iv)*] this solution is unique in* $C^2_{\mathrm{pol}}$ *(resp.* $C^2_{\mathrm{pol},+}$; $C^2_b$*). Moreover, there exists a Markov control policy that is optimal for all* $x \in \mathbb{R}^I$.

We end this section with a few remarks on the role of the treelike assumption in this work. Our results strongly depend on estimates on moments of the controlled process that are subexponential in the time variable [in cases (i) and (ii) of the main result]. These are obtained by considering a deterministic model in place of (18)–(22). Use $w$ in place of $x + \widetilde{W}$ [where as in (18) $x$ is the initial condition for $X$] and use $x = x(t)$ (resp. $y$, $z$, $\psi$) in place of $X$ (resp. $Y$, $Z$, $\Psi$). Then

$$(33) \qquad x_i(t) = w_i(t) - \sum_{j \in \mathcal{J}} \mu_{ij} \int_0^t \psi_{ij}(s)\,ds - \theta_i \int_0^t y_i(s)\,ds, \qquad i \in \mathcal{I},$$

$$(34) \qquad \sum_{j \in \mathcal{J}} \psi_{ij} = x_i - y_i, \qquad i \in \mathcal{I},$$

$$(35) \qquad \sum_{i \in \mathcal{I}} \psi_{ij} = -z_j, \qquad j \in \mathcal{J},$$

(36)     $\quad y_i, z_j \geq 0, \qquad i \in \mathcal{I}, j \in \mathcal{J},$

(37)     $\quad e \cdot y \wedge e \cdot z = 0,$

where $\psi_{ij} = 0$, $i \not\sim j$. The first two cases of Theorem 1 are based on showing that (33)–(37) imply an estimate of the form

(38)     $$\|x(t)\| \leq m_0 (1+t)^{m_0} (1 + \|w\|_t^*)^{m_0},$$

where $m_0$ does not depend on $t$, $w$, $y$, $z$ and $\psi$. If the treelike assumption is removed, this estimate does not hold true in general (as shown in the example below). We leave open the question of whether this implication holds true in full generality under the treelike assumption. We stress that this is the reason for imposing the treelike assumption in this paper, rather than the using capability to rephrase the model equations (18)–(22) as (23): Although (23) is useful, its absence would not be a serious obstacle to treating the problem, whereas the large time estimates constitute a key ingredient of the proof.

EXAMPLE 1.   Consider a system with classes 1 and 2 and stations $A$ and $B$, and with $\mu_{1A} = \mu_{2A} = 1$, $\mu_{1B} = \mu_{2B} = 2$, and arbitrary $\theta_1$ and $\theta_2$. Consider $w = 0$, $\psi_{1A} = -\psi_{2A} = k$, $\psi_{1B} = -\psi_{2B} = -k(1 + e^{-2t})/2$, $x_1 = -x_2 = k(1 - e^{-2t})/2$, $y = 0$ and $z = 0$. Then (33)–(37) hold for every $k > 0$. Thus (38) cannot hold.

Finally, there is another central role played by the treelike condition, as elaborated in [2]. The diffusion model turns out to depend on whether preemption is allowed in the queueing system under scaling. As explained in [2], preemptive and nonpreemptive policies give rise to genuinely different diffusion models if the treelike assumption does not hold, whereas under the treelike assumption the corresponding diffusion models coincide (as supported by the result in [4] for the case of a single station and in [2] for cases (i) and (ii) of Theorem 1).

**3. Estimating the state $X$ in terms of $e \cdot X$.**  While the relationship $\|Y(t)\| + \|Z(t)\| \leq c\|X(t)\|$ is immediate from (25) and (26), the following result shows that in a weaker sense $Y$, $Z$ and $x + \widetilde{W}$ dominate $X$ (or in the deterministic notation, $y$, $z$ and $w$ dominate $x$). The result uses only the relationships (33)–(35), and not the further constraints (36) and (37).

PROPOSITION 1.   *Let* (33)–(35) *hold. Then there is a constant $m_1$ not depending on $\psi, w, x, y, z$ or $t$, such that*

$$\|\mathfrak{I}\psi(t)\| + \|x(t)\| \leq m_1 (\|w\|_t^* + \|\mathfrak{I}y\|_t^* + \|\mathfrak{I}z\|_t^*), \qquad t \geq 0.$$

Note that if, in addition, (36) and (37) are assumed, then $\|\mathfrak{I}y\|_t^* = \mathfrak{I}(e \cdot x)^+(t)$ and $\|\mathfrak{I}z\|_t^* = \mathfrak{I}(e \cdot x)^-(t)$. As a result of Proposition 1, the state $x$ is dominated by $w$ and $e \cdot x$ in the sense

$$\|x(t)\| \le m_1(\|w\|_t^* + \mathfrak{I}|e \cdot x|(t)).$$

LEMMA 1. *Let $w$ be a measurable, locally bounded function and assume $x = w - \mu \int_0^t x(s)\,ds$. Then*

$$(39) \qquad x(t) = w(t) - \mu \int_0^t w(s)e^{-\mu(t-s)}\,ds$$

*and, in particular, if $\mu > 0$, then $|x(t)| \le 2|w|_t^*$, $t \ge 0$.*

PROOF. Uniqueness of solutions $x$ is standard and (39) is verified by substitution. $\square$

PROOF OF PROPOSITION 1. Observe, by replacing $w_i - \theta_i \mathfrak{I}y_i$ with $w_i$, that without loss we can take $\theta_i = 0$, $i \in \mathcal{I}$; therefore, $\theta_i = 0$ in the sequel.

A node in a tree is said to be a leaf if there is exactly one edge joining it. Recall that the tree $\mathcal{T}$ has $\kappa = I + J$ nodes and $\kappa - 1$ edges. Let $\mathcal{T}_1, \mathcal{T}_2, \ldots, \mathcal{T}_{\kappa-1} = \mathcal{T}$ be an increasing sequence of trees as follows. For $n = 1, \ldots, k-2$, $\mathcal{T}_n$ is obtained from $\mathcal{T}_{n+1}$ by deleting a leaf and the edge joining it. Note that $\mathcal{T}_1$ has exactly two nodes: one in $\mathcal{I}$ and one in $\mathcal{J}$. Let $\mathcal{V}_n$ denote the vertex set of $\mathcal{T}_n$. Let $v_{n+1} = \mathcal{V}_{n+1} \setminus \mathcal{V}_n$ denote the node in $\mathcal{V}_{n+1}$ that does not belong to $\mathcal{V}_n$.

Denote $\mathcal{I}_n = \mathcal{I} \cap \mathcal{V}_n$ and $\mathcal{J}_n = \mathcal{J} \cap \mathcal{V}_n$. We show that for $n = 1, \ldots, \kappa - 1$, if

$$x_i = w_i - \sum_{j \in \mathcal{J}_n} \mu_{ij} \mathfrak{I}\psi_{ij}, \qquad i \in \mathcal{I}_n,$$

$$(40) \qquad \sum_{j \in \mathcal{J}_n} \psi_{ij} = x_i + \alpha_i, \qquad i \in \mathcal{I}_n,$$

$$\sum_{i \in \mathcal{I}_n} \psi_{ij} = \beta_j, \qquad j \in \mathcal{J}_n,$$

then

$$(41) \quad \sum_{\substack{i \in \mathcal{I}_n \\ j \in \mathcal{J}_n}} |\mathfrak{I}\psi_{ij}|_t^* + \sum_{i \in \mathcal{I}_n} |x_i|_t^* \le c_n\left( \sum_{i \in \mathcal{I}_n} (|w_i|_t^* + |\mathfrak{I}\alpha_i|_t^*) + \sum_{j \in \mathcal{J}_n} |\mathfrak{I}\beta_j|_t^* \right).$$

The implication $(40) \Rightarrow (41)$ is proved by induction on $n$.

*Induction base*: $n = 1$. $\mathcal{T}_1$ has exactly two nodes, say $i \in \mathcal{I}$ and $j \in \mathcal{J}$. By (40), $x_i = w_i - \mu_{ij}\Im\psi_{ij}$ and $\Im\psi_{ij} = \Im\beta_j$. Hence (41) holds.

*Induction step*. Assume that $(40) \Rightarrow (41)$ holds for $n \in [1, \kappa - 2]$. Let (40) hold for $n + 1$. We show that (41) holds for $n + 1$ in the following two cases.

CASE 1.    *The leaf node $v_{n+1}$ is in $\mathcal{I}$*. Denote $i_0 = v_{n+1}$ and let $j_0$ denote the unique node $j \sim i_0$ in $\mathcal{T}_{n+1}$. The validity of (40) for $n + 1$ implies

$$
\begin{aligned}
x_i &= w_i - \sum_{j \in \mathcal{J}_n} \mu_{ij}\Im\psi_{ij}, \qquad i \in \mathcal{I}_n, \\
\sum_{j \in \mathcal{J}_n} \psi_{ij} &= x_i + \alpha_i, \qquad i \in \mathcal{I}_n, \\
\sum_{i \in \mathcal{I}_n} \psi_{ij} &= \beta_j, \qquad j \in \mathcal{J}_n \setminus \{j_0\}, \\
\sum_{i \in \mathcal{I}_n} \psi_{ij_0} &= \beta_{j_0} - \psi_{i_0 j_0},
\end{aligned}
\tag{42}
$$

regarding $i \in \mathcal{I}_n$, $j \in \mathcal{J}_n$, and

$$
x_{i_0} = w_{i_0} - \mu_{i_0 j_0}\Im\psi_{i_0 j_0}, \qquad \psi_{i_0 j_0} = x_{i_0} + \alpha_{i_0}.
\tag{43}
$$

regarding the node $i_0$. By (42) and the induction assumption,

$$
\begin{aligned}
&\sum_{i \in \mathcal{I}_n, j \in \mathcal{J}_n} |\Im\psi_{ij}|_t^* + \sum_{i \in \mathcal{I}_n} |x_i|_t^* \\
&\qquad \leq c_n \left( \sum_{i \in \mathcal{I}_n} (|w_i|_t^* + |\Im\alpha_i|_t^*) + \sum_{j \in \mathcal{J}_n} |\Im\beta_j|_t^* + |\Im\psi_{i_0 j_0}|_t^* \right).
\end{aligned}
\tag{44}
$$

By (43),

$$
x_{i_0} = w_{i_0} - \mu_{i_0 j_0}\Im\alpha_{i_0} - \mu_{i_0 j_0}\Im x_{i_0}.
\tag{45}
$$

Applying Lemma 1 to (45) and again using (43) shows

$$
|x_{i_0}(t)| + |\Im\psi_{i_0 j_0}(t)| \leq c'(|w_{i_0}|_t^* + |\Im\alpha_{i_0}|_t^*).
\tag{46}
$$

Combining (44) and (46) establishes the validity of (41) for $n + 1$.

CASE 2. *The leaf node $v_{n+1}$ is in $\mathcal{J}$.* Denote $j_0 = v_{n+1}$ and let $i_0$ denote the unique node $i \sim j_0$ in $\mathcal{T}_{n+1}$. Assuming (40) for $n+1$ implies

$$
\begin{aligned}
x_i &= w_i - \sum_{j \in \mathcal{J}_n} \mu_{ij} \mathfrak{I} \psi_{ij}, & i \in \mathcal{I}_n \setminus \{i_0\}, \\
\sum_{j \in \mathcal{J}_n} \psi_{ij} &= x_i + \alpha_i, & i \in \mathcal{I}_n \setminus \{i_0\}, \\
x_{i_0} &= (w_{i_0} - \mu_{i_0 j_0} \mathfrak{I} \psi_{i_0 j_0}) - \sum_{j \in \mathcal{J}_n} \mu_{i_0 j} \mathfrak{I} \psi_{i_0 j}, \\
\sum_{j \in \mathcal{J}_n} \psi_{i_0 j} &= x_{i_0} + (\alpha_{i_0} - \psi_{i_0 j_0}), \\
\sum_{i \in \mathcal{I}_n} \psi_{ij} &= \beta_j, & j \in \mathcal{J}_n,
\end{aligned}
\tag{47}
$$

and

$$
\psi_{i_0 j_0} = \beta_{j_0}.
\tag{48}
$$

By (47) and the induction assumption, (44) follows. Combining (44) with (48) gives (41) for $n+1$.

This completes the proof that (40) implies (41) for $n \in [1, \kappa - 1]$. The result follows on taking $n = \kappa - 1$ and substituting $\alpha = -y$ and $\beta = -z$. □

## 4. An integral formula for $Y$ and $Z$.

Equations (33)–(35) were used in the previous section as substitutes for (18)–(20). They express a relationship between the quantities $x$, $y$, $z$, $w$ and $\psi$. In this section we extract a relationship between $y$, $z$ and $w$ alone. This relationship, in the form of an integral equation, is a key element in treating Theorem 1(i) and (ii).

For $\alpha \in \mathbb{R}$ denote

$$
\mathfrak{T}_\alpha f = f + \alpha \mathfrak{I} f.
$$

Lemma 1 shows that $\mathfrak{T}_\alpha$ is invertible. It is easy to see that the operators $\mathfrak{T}_\alpha$ and $\mathfrak{T}_\beta$ commute. If $A = (\alpha_1, \dots, \alpha_k)$ is a finite real-valued sequence, denote

$$
\mathfrak{T}_A = \mathfrak{T}_{\alpha_1} \circ \cdots \circ \mathfrak{T}_{\alpha_k}.
$$

Then $\mathfrak{T}_A$ does not depend on the order of the elements of $A$, but it depends on the multiplicity of each element. Set $\mathfrak{T}_\varnothing$ corresponding to $k = 0$ as the identity map. Equations (33)–(35) imply

$$
\sum_{j \in \mathcal{J}} \mathfrak{T}_{\mu_{ij}} \psi_{ij} = w_i - \mathfrak{T}_{\theta_i} y_i, \qquad i \in \mathcal{I},
\tag{49}
$$

$$
\sum_{i \in \mathcal{I}} \psi_{ij} = -z_j, \qquad j \in \mathcal{J}.
\tag{50}
$$

THEOREM 2. *Let* (49) *and* (50) *hold.*

(i) *The quantities $y$ and $z$ solve the integral equation*

$$\sum_{i \in \mathcal{I}} \mathfrak{T}_{A_i} w_i - \sum_{i \in \mathcal{I}} \mathfrak{T}_{A_i'} y_i + \sum_{j \in \mathcal{J}} \mathfrak{T}_{B_j} z_j = 0, \tag{51}$$

*where $A_i$ and $B_j$ are finite (possibly empty) sequences with values in $\{\mu_{ij} : (i, j) \in \mathcal{E}\}$ and, for $i \in \mathcal{I}$, $A_i'$ is the concatenation of $A_i$ with $\theta_i$.*

(ii) *In the special case where $\mu_{ij} = \mu_j$ for $(i, j) \in \mathcal{E}$ and $\theta_i = 0$, $i \in \mathcal{I}$, equation* (51) *takes the form*

$$\sum_{i \in \mathcal{I}} (w_i - y_i) + \sum_{j \in \mathcal{J}} \mathfrak{T}_{\mu_j} z_j = 0. \tag{52}$$

(iii) *In the special case where $\mu_{ij} = \mu_i$ for $(i, j) \in \mathcal{E}$ and $\theta_i = 0$, $i \in \mathcal{I}$, equation* (51) *takes the form*

$$\sum_{i \in \mathcal{I}} \mathfrak{T}_{M_i} (w_i - y_i) + \mathfrak{T}_M (e \cdot z) = 0, \tag{53}$$

*where $M_i = (\mu_{i'})_{i' \in \mathcal{I}, \ i' \neq i}$ and $M = (\mu_{i'})_{i' \in \mathcal{I}}$.*

REMARK 1. (a) Writing $\mathfrak{T}_n$ for the $n$-power $\mathfrak{T} \circ \cdots \circ \mathfrak{T}$ of the operator $\mathfrak{T}$, it is useful to note that the integral equation (51) can be written as

$$\begin{aligned}
e \cdot w - e \cdot y + e \cdot z \\
+ \sum_{i \in \mathcal{I}} \sum_{n=1}^{m_i} a_{i,n} \mathfrak{T}_n w_i - \sum_{i \in \mathcal{I}} \sum_{n=1}^{m_i'} a_{i,n}' \mathfrak{T}_n y_i + \sum_{j \in \mathcal{J}} \sum_{n=1}^{m_j''} a_{j,n}'' \mathfrak{T}_n z_j = 0.
\end{aligned} \tag{54}$$

Here, $m_i, m_i', m_j''$, $a_{i,n}$, $a_{i,n}'$ and $a_{j,n}''$ are positive constants that we do not give in explicit form.

(b) Recall that under (36) and (37) we have $y = (e \cdot x)^+ u$ and $z = (e \cdot x)^- v$. As a result, (54) expresses a relationship between the data $w$, the controls $u$ and $v$, and the quantity $e \cdot x$ alone.

We need some notation regarding the tree $\mathcal{T}$ to be used in this and the following sections. Fix one of the class nodes, $i_0$, as a root. (Analogous notation applies if we fix a station node, some $j_0$ as a root.) For $k = 0, 1, \ldots$, let level $k$, denoted by $l_k$, be the set of nodes of $\mathcal{T}$ at distance $k$ from the root $i_0$ along the edges of $\mathcal{T}$ (see Figure 2). Note that $l_0 = \{i_0\}$ and that $l_k$ is empty for all $k$ large. Let also

$$L_k = l_0 \cup l_1 \cup \cdots \cup l_k$$

be the set of nodes at distance at most $k$ from the root and let

$$L_k^{\mathcal{I}} = L_k \cap \mathcal{I}, \qquad L_k^{\mathcal{J}} = L_k \cap \mathcal{J}.$$

Note that the elements of $L_k^{\mathcal{I}}$ (resp. $L_k^{\mathcal{J}}$) are at even (resp. odd) distance from the root, not exceeding $k$. Let $K$ be the largest $k$ such that $l_k$ is nonempty. For a node $v$ at level $k$ let $B(v)$ ($B$ for below) be the set of nodes $v' \sim v$ at level $k+1$. For a node $v$ at level $k \in [1, K]$ let $a(v)$ ($a$ for above) be the unique node $v' \sim v$ at level $k-1$.

PROOF OF THEOREM 2. (i) Suppose we prove (51) for the case $\theta_i = 0$, $i \in \mathcal{I}$, that is,

$$(55) \qquad \sum_{i \in \mathcal{I}} \mathfrak{T}_{A_i}(w_i - y_i) + \sum_{j \in \mathcal{J}} \mathfrak{T}_{B_j} z_j = 0.$$

Then, for arbitrary $\theta_i$, (51) is obtained from (55) on substituting $\mathfrak{T}_{\theta_i} y_i$ for $y_i$. Therefore, in the sequel we set $\theta_i = 0$, $i \in \mathcal{I}$, and turn to prove (55).

We show that for $k \geq 1$, we have

$$(56) \qquad \sum_{i \in L_{2k-2}^{\mathcal{I}}} \mathfrak{T}_{A_i^k}(w_i - y_i) + \sum_{j \in L_{2k-1}^{\mathcal{J}}} \mathfrak{T}_{B_j^k} z_j + \sum_{i \in l_{2k}} \mathfrak{T}_{C_i^k} \psi_{ia(i)} = 0,$$

where $A_i^k$, $B_j^k$ and $C_i^k$ are (possibly empty) sequences with values in $\{\mu_{ij} : (i,j) \in \mathcal{E}\}$, and summation over an empty set is regarded as zero. Equation (55) follows on taking $k$ larger than $K$.

We prove (56) by induction on $k$.

*Induction base*: $k = 1$. By (49),

$$\sum_{j \in l_1} \mathfrak{T}_{\mu_{i_0 j}} \psi_{i_0 j} = w_{i_0} - y_{i_0}$$

and by (50),

$$\psi_{i_0 j} = -z_j - \sum_{i \in l_2} \psi_{ij}, \qquad j \in l_1.$$

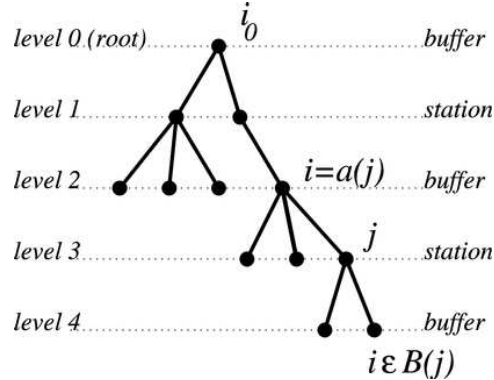

FIG. 2. *Buffer-station tree.*

It follows that

$$w_{i_0} - y_{i_0} + \sum_{j \in l_1} \mathfrak{T}_{\mu_{i_0 j}} z_j + \sum_{i \in l_2} \mathfrak{T}_{\mu_{i_0 a(i)}} \psi_{ia(i)} = 0$$

and (56) holds.

*Induction step.* Assume that (56) holds for $k$. Using (49) and then (50), for $i \in l_{2k}$ and $j = a(i)$ we have

$$\mathfrak{T}_{\mu_{ij}} \psi_{ij} = -\sum_{j' \in B(i)} \mathfrak{T}_{\mu_{ij'}} \psi_{ij'} + w_i - y_i$$

$$= \sum_{j' \in l_{2k+1}} \left\{ \mathfrak{T}_{\mu_{ij'}} z_{j'} + \sum_{i' \in l_{2k+2}} \mathfrak{T}_{\mu_{ij'}} \psi_{i'j'} \right\} + w_i - y_i.$$

Apply $\mathfrak{T}_{C_i^k}$ on the above equation [where still $i \in l_{2k}$ and $j = a(i)$] to get

$$
\begin{aligned}
(57) \qquad \mathfrak{T}_{C_i^k} \mathfrak{T}_{\mu_{ij}} \psi_{ij} &= \sum_{j' \in l_{2k+1}} \mathfrak{T}_{C_i^k} \mathfrak{T}_{\mu_{ij'}} z_{j'} \\
&\quad + \sum_{i' \in l_{2k+2}} \mathfrak{T}_{C_i^k} \mathfrak{T}_{\mu_{ia(i')}} \psi_{i'a(i')} + \mathfrak{T}_{C_i^k}(w_i - y_i).
\end{aligned}
$$

Let $D_k = (\mu_{ia(i)})_{i \in l_{2k}}$. Apply $\mathfrak{T}_{D_k}$ to (56) and use (57) to substitute for each summand in the third sum in (56) to get

$$
\begin{aligned}
0 &= \sum_{i \in L_{2k-2}^{\mathcal{I}}} \mathfrak{T}_{D_k} \mathfrak{T}_{A_i^k}(w_i - y_i) + \sum_{j \in L_{2k-1}^{\mathcal{J}}} \mathfrak{T}_{D_k} \mathfrak{T}_{B_j^k} z_j \\
&\quad + \sum_{i \in l_{2k}} \left\{ \sum_{j' \in l_{2k+1}} \mathfrak{T}_{D_k \setminus \{\mu_{ij'}\}} \mathfrak{T}_{\mu_{ij'}} \mathfrak{T}_{C_i^k} z_{j'} \right. \\
&\quad \left. + \sum_{i' \in l_{2k+2}} \mathfrak{T}_{D_k \setminus \{\mu_{ij}\}} \mathfrak{T}_{\mu_{ia(i')}} \mathfrak{T}_{C_i^k} \psi_{i'a(i')} + \mathfrak{T}_{D_k \setminus \{\mu_{ij}\}} \mathfrak{T}_{C_i^k}(w_i - y_i) \right\},
\end{aligned}
$$

where $D_k \setminus \{\mu\}$ denotes a sequence obtained from $D_k$ by deleting from it one element of value $\mu$. This proves that (56) holds for $k + 1$ and completes the proof of part (i).

(ii) In the case that $\mu_{ij} = \mu_j$ and $\theta_i = 0$, by (49) and (50),

$$\sum_{j \in \mathcal{J}} \mathfrak{T}_{\mu_j} \psi_{ij} = w_i - y_i, \qquad i \in \mathcal{I},$$

$$\sum_{i \in \mathcal{I}} \mathfrak{T}_{\mu_j} \psi_{ij} = -\mathfrak{T}_{\mu_j} z_j, \qquad j \in \mathcal{J},$$

hence

$$\sum_{i \in \mathcal{I}} (w_i - y_i) + \sum_{j \in \mathcal{J}} \mathfrak{T}_{\mu_j} z_j = 0.$$

(iii) In the case that $\mu_{ij} = \mu_i$ and $\theta_i = 0$, by (49) and (50),

$$(58) \qquad \mathfrak{T}_{\mu_i} \sum_{j \in \mathcal{J}} \psi_{ij} = w_i - y_i, \qquad i \in \mathcal{I},$$

$$(59) \qquad \sum_{i \in \mathcal{I}} \psi_{ij} = -z_j, \qquad j \in \mathcal{J}.$$

The result follows on applying $\prod_{i' \neq i} \mathfrak{T}_{\mu_{i'}}$ on (58) and $\prod_{i' \in \mathcal{I}} \mathfrak{T}_{\mu_{i'}}$ on (59). $\quad\square$

PROPOSITION 2. *Let* (49) *and* (50) *hold. Assume also that*

$$(60) \qquad y_i \geq 0, \qquad z_j \geq 0, \qquad i \in \mathcal{I}, j \in \mathcal{J},$$

$$(61) \qquad e \cdot y(t) \wedge e \cdot z(t) = 0.$$

*In cases* (ii) *and* (iii) *of Theorem* 2 *there is a constant $m_2$ such that*

$$(62) \qquad \|y(t)\| + \|z(t)\| \leq m_2 (1+t)^{m_2} \|w\|_t^*, \qquad t \geq 0.$$

PROOF. In the case $\mu_{ij} = \mu_j$, (52) can be written as

$$e \cdot w - e \cdot y + e \cdot z + \sum_{j \in \mathcal{J}} \mu_j \mathfrak{I} z_j = 0.$$

By (60), $z_j$ are positive and therefore $e \cdot w - e \cdot y + e \cdot z \leq 0$. Thus by (61), $0 \leq e \cdot z \leq (-e \cdot w)^+$ and, therefore,

$$e \cdot w \leq e \cdot y - e \cdot z \leq e \cdot w + t |e \cdot w|_t^*.$$

Since by (61), $\|y\| + \|z\| = |e \cdot y - e \cdot z|$, (62) follows.

In the case $\mu_{ij} = \mu_i$, applying $\mathfrak{T}_{\mu_i}^{-1}$ to (58) and by (59) [or by applying $\mathfrak{T}_M^{-1}$ to (53)], we have

$$\sum_{i \in \mathcal{I}} \mathfrak{T}_{\mu_i}^{-1}(w_i - y_i) + e \cdot z = 0;$$

hence, by Lemma 1 and positivity of $y_i$,

$$e \cdot w - e \cdot y + e \cdot z = \sum_{i \in \mathcal{I}} \mu_i \int_0^t (w_i(s) - y_i(s)) \exp(-\mu_i(t-s)) \, ds \leq c_1 t \|w\|_t^*$$

for some constant $c_1$. By positivity of $z_j$ and (61) we, therefore, have

$$(63) \qquad \|z(t)\| \leq c_2 (1+t) \|w\|_t^*.$$

By (54), the positivity of its coefficients and of $y_i$, and (63), there is a constant $c_3$ such that

$$(64) \qquad \|y(t)\| \leq c_3 (1+t)^{c_3} \|w\|_t^*.$$

Combining (63) and (64) establishes (62). $\quad\square$

COROLLARY 1. *Under the assumptions of Theorem* 1(i), *for any* $m \geq 1$, *any initial condition* $x \in \mathbb{R}^I$ *and any admissible system* $\pi \in \Pi$,

$$E_x^\pi \|X(t)\|^m \leq m_3 (1 + \|x\|)^{m_3} (1 + t)^{m_3}, \qquad t \geq 0,$$

*where* $m_3$ *do not depend on* $x$, $\pi$ *and* $t$.

PROOF. By Propositions 1 and 2,

$$\|x(t)\| \leq m_1 (\|w\|_t^* + \|\Im y\|_t^* + \|\Im z\|_t^*)$$
$$\leq [m_1 + m_1 m_2 t (1 + t)^{m_2}] \|w\|_t^*.$$

Thus $E_x^\pi \|X(t)\|^m \leq c(1 + t)^c E_x^\pi (\|x\| + c \|\widetilde{W}\|_t^*)^m$ for a constant $c$ depending only on $m_1$, $m_2$ and $m$, and the result follows from standard estimates on the Brownian motion.  □

## 5. The nonidling property.

This section investigates a relationship between a property defined below and referred to as the nonidling property, and the uniform estimate

$$\|x(t)\| \leq m_4 (1 + t)^{m_4} (1 + \|w\|_t^*), \tag{65}$$

where $m_4$ does not depend on $t$, $w$, $y$, $z$ and $\psi$. In particular, using the integral equation developed in the previous section, it shows that for trees of diameter not exceeding 3, this property implies the uniform estimate. A relatively simple argument then shows that the property holds for such trees and the estimate follows.

Rewrite relationships (33)–(35) as in the previous section and recall relationships (36) and (37):

$$\sum_{j \in \mathcal{J}} (\psi_{ij} + \mu_{ij} \Im \psi_{ij}) = w_i - y_i - \theta_i \Im y_i, \qquad i \in \mathcal{I}, \tag{66}$$

$$\sum_{i \in \mathcal{I}} \psi_{ij} = -z_j, \qquad j \in \mathcal{J}, \tag{67}$$

$$y_i, z_j \geq 0, \qquad i \in \mathcal{I}, j \in \mathcal{J}, \tag{68}$$

$$e \cdot y \wedge e \cdot z = 0. \tag{69}$$

Note that $x$ [cf. (33)] is not a part of this system of equations, but can be obtained from it via

$$x_i = w_i - \sum_{j \in \mathcal{J}} \mu_{ij} \Im \psi_{ij} - \theta_i \Im y_i.$$

We say that *the system* (66)–(69) *incurs no idleness on* $[0, T)$ (cf. [18]) if $z_j(t) = 0$ for $t \in [0, T)$, $j \in \mathcal{J}$.

THE NONIDLING PROPERTY. For every $T > 0$, if the system starts with $w_i(0) > 0$, and $w_i$, $i \in \mathcal{I}$, are strictly increasing and right-continuous on $[0, T)$, then the system incurs no idleness on $[0, T)$.

A tree $\mathcal{T}$ of diameter 3 has the form depicted in Figure 3, where there are only two nodes $i_0 \in \mathcal{I}$ and $j_0 \in \mathcal{J}$ that are not leaves.

THEOREM 3. *Let Assumption 1 hold and assume the diameter of the tree $\mathcal{T}$ is 3. Then the nonidling property [for the system (66)–(69)] implies the estimate (65), where $m_4$ does not depend on $t$, $w$, $y$, $z$ and $\psi$.*

LEMMA 2. *Let the assumptions of Theorem 3 hold. Fix $T > 0$ and let $(\psi, y, z, w)$ be given, satisfying (66)–(69) on $[0, T]$. Let bounded measurable functions $f_i \geq 0$, $i \in \mathcal{I}$, be given. Then there is a constant $m_5$ that does not depend on $T$, $f_i$, $\psi$, $y$, $z$ or $w$, and there exist $(\hat{\psi}, \hat{y}, \hat{z}, \hat{w})$ defined on $[0, T]$ that satisfy (66)–(69), and, moreover,*

$$(70) \qquad\qquad e \cdot \hat{y} \geq e \cdot y, \qquad e \cdot \hat{z} \leq e \cdot z$$

*and*

$$(71) \qquad\qquad \hat{w}_i = w_i + f_i + \eta_i,$$

*where $\eta_i$ are nondecreasing and continuous, and $0 \leq \eta_i \leq m_5 T \sum_{i' \in \mathcal{I}} |f_{i'}|_T^*$ on $[0, T]$.*

PROOF. It suffices to prove the lemma in the case where all but one of the functions $f_i$ vanish, since the argument can then be repeated. Consider first the case where $f_i$ vanish for all $i \neq i_0$ (and $f_{i_0} \geq 0$). Define $\hat{w}_i = w_i$, $\hat{y}_i = y_i$ and $\hat{\psi}_{ij} = \psi_{ij}$ for $i \neq i_0$. Define now $\hat{w}_{i_0 j}$ and $\hat{z}_j$ as follows. For $t \in \Theta := \{s \in [0, T] : e \cdot z(s) \geq f_{i_0}(s)\}$, let $j = I + 1$, let $\Delta_j := z_j \wedge f_{i_0}$, and set $\hat{\psi}_{i_0 j} = \psi_{i_0 j} + \Delta_j$ and $\hat{z}_j = z_j - \Delta_j$. Similarly, for $j \in [I + 2, I + J]$, let $\Delta_j = z_j \wedge (f_{i_0} - \Delta_{I+1} - \cdots - \Delta_{j-1})$ and $\hat{\psi}_{i_0 j} = \psi_{i_0 j} + \Delta_j$, $\hat{z}_j = z_j - \Delta_j$. Note that for all $j$, $\hat{z}_j \geq 0$ and that (67) is satisfied by the hat system. It is elementary to check that

$$(72) \qquad \Delta_j(t) \geq 0, \qquad j \in \mathcal{J}; \qquad \sum_{j \in \mathcal{J}} \Delta_j(t) = f_{i_0}, \qquad t \in \Theta.$$

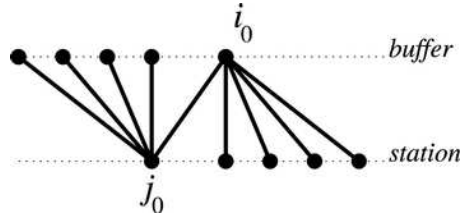

FIG. 3. *Tree of diameter 3.*

As a result, with $\delta(t) = 0$,

$$(73) \qquad \sum_{j \in \mathcal{J}} \hat{\psi}_{i_0 j}(t) = \sum_{j \in \mathcal{J}} \psi_{i_0 j}(t) + f_{i_0}(t) - \delta(t).$$

For $t \in \Theta^c$, set $\hat{z}_j = 0$ and $\hat{\psi}_{i_0 j} = \psi_{i_0 j} + \Delta_j$, where now $\Delta_j = z_j$, $j \in \mathcal{J}$. Then again (67) holds for the hat system. Also (73) holds with $\delta(t) = f_{i_0}(t) - e \cdot z(t)$. Note that

$$(74) \qquad 0 \leq \delta \leq f_{i_0}.$$

Now that $\hat{\psi}$ is defined on $[0, T]$, let

$$\hat{w}_{i_0} = w_{i_0} + \sum_{j \in \mathcal{J}} (\hat{\psi}_{i_0 j} - \psi_{i_0 j} + \mu_{i_0 j} \mathfrak{I}(\hat{\psi}_{i_0 j} - \psi_{i_0 j})) + \delta + \theta_{i_0} \mathfrak{I}\delta$$

and let $\hat{y}_{i_0} = y_{i_0} + \delta$. Then (66) holds on $[0, T]$ for $i = i_0$. By (73), $\hat{w}_{i_0} = w_{i_0} + f_{i_0} + \eta_{i_0}$, where $\eta_{i_0} = \sum_{j \in \mathcal{J}} \mu_{i_0 j} \mathfrak{I}(\hat{\psi}_{i_0 j} - \psi_{i_0 j}) + \theta_{i_0} \mathfrak{I}\delta$. By (72), $0 \leq \hat{\psi}_{i_0 j} - \psi_{i_0 j} \leq f_{i_0}$ on $\Theta$ and by construction the same holds on $\Theta^c$. Combined with (74), the properties of $\eta_{i_0}$ stated in the lemma follow. Inequalities (70) also hold by construction. Condition (69) holds for the hat system since on $\Theta$, $\hat{y} = y$ and by (72) $e \cdot \hat{z} \leq e \cdot z$; and on $\Theta^c$, $e \cdot \hat{z} = 0$.

Next consider the case where for some $i_1 \neq i_0$, $f_i = 0$ on $[0, T]$ for all $i \neq i_1$. The argument is similar, but slightly more complicated because $i_1 \not\sim j$ for most $j$. We indicate only where the argument changes. The definition of $\hat{\psi}_{ij}$ and $\hat{z}_j$ differs as follows. Variables $\Delta_j$, $\delta$ and $\hat{z}_j$ are defined as before, where $i_0$ is replaced by $i_1$. Recall that in the previous case, $\hat{\psi}$ is defined as $\hat{\psi} = \psi + \Delta$. In the current case, let $\hat{\psi}_{i_1 j_0} = \psi_{i_1 j_0} + \sum_j \Delta_j$, $\hat{\psi}_{i_0 j} = \psi_{i_0 j} + \Delta_j$ for $j \neq j_0$, $\hat{\psi}_{i_0 j_0} = \psi_{i_0 j_0} - \sum_{j \neq j_0} \Delta_j$ and $\hat{\psi}_{ij} = \psi_{ij}$ for $i \notin \{i_0, i_1\}$. Set

$$\hat{w}_{i_1} = w_{i_1} + \hat{\psi}_{i_1 j_0} - \psi_{i_1 j_0} + \mu_{i_1 j_0} \mathfrak{I}(\hat{\psi}_{i_1 j_0} - \psi_{i_1 j_0}) + \delta + \theta_{i_1} \mathfrak{I}\delta,$$

$$\hat{w}_{i_0} = w_{i_0} + \sum_{j \neq j_0} \mu_{i_0 j} \mathfrak{I}(\hat{\psi}_{i_0 j} - \psi_{i_0 j}) =: w_{i_0} + \tilde{\eta}_{i_0}.$$

Let $\hat{y}_{i_1} = y_{i_1} + \delta$ and let $\hat{y}_i = y_i$ for all $i \neq i_1$. We check that all properties (66)–(69) are satisfied by the hat system and that the conclusions of the lemma all hold, except that $\hat{w}_{i_0} = w_{i_0} + f_{i_0} + \tilde{\eta}_{i_0}$, where there is no guarantee that $\tilde{\eta}_{i_0}$ is nondecreasing (in fact, we have assumed without loss that $f_{i_0} = 0$). However, this is now corrected by applying a further modification as in the previous paragraph, where now we take $f_{i_0} = \sum_{j \neq j_0} \mu_{i_0 j} \Im[(\hat{\psi}_{i_0 j} - \psi_{i_0 j})^-]$, so that overall $f_{i_0} + \tilde{\eta}_{i_0}$ is nondecreasing and $f_i = 0$ for all $i \neq i_0$. $\square$

PROOF OF THEOREM 3. Fix $T > 0$ and let $(\psi, y, z, w)$ be given, satisfying (66)–(69) on $[0, T+1]$. Let $f_i$ be defined on $[0, T+1]$ as

$$f_i(t) = 1 + \max\left(0, \sup_{0 \leq s < T+1} w_i(s)\right) - w_i(t) + t/T =: \bar{a}_i(w, T) - w_i(t) + t/T.$$

Let $(\hat{\psi}, \hat{y}, \hat{z}, \hat{w})$ be as in the conclusion of the lemma. Then $\hat{w}_i = \bar{a}_i(w, T) + t/T + \eta_i$; thus $\hat{w}_i(0) > 0$ and $\hat{w}_i$ is strictly increasing on $[0, T+1)$. Therefore, by the nonidling property, $\hat{z}_j = 0$ on $[0, T+1)$, $j \in \mathcal{J}$. Applying (54) to $\hat{w}, \hat{y}, \hat{z}$, and using the nonnegativity of $y_i$ and of the constants $a'_{i,n}$,

$$\sum_{i \in \mathcal{I}} \hat{y}_i(t) \leq c_1 (1+t)^{c_1} \sum_{i \in \mathcal{I}} |\hat{w}_i|_T^*, \qquad t \leq T,$$

where $c_1$ does not depend on $\hat{w}$ and $T$. Using (70) and (71),

$$e \cdot y(t) \leq c_2 (1+t)^{c_2} \sum_{i \in \mathcal{I}} (\bar{a}_i(w, T) + \eta_i(T)) \leq c_3 (1+t)^{c_3} T (1 + \|w\|_T^*), \qquad t \leq T.$$

Thus on $[0, T]$,

$$e \cdot y \leq c_4 (1+T)^{c_4} (1 + \|w\|_T^*).$$

Using again (54), now on $w, y, z$ and equipped with the bound on $y$,

$$e \cdot z \leq c_5 (1+T)^{c_5} (1 + \|w\|_T^*).$$

As a result,

$$\|y\|_T^* + \|z\|_T^* \leq c_6 (1+T)^{c_6} (1 + \|w\|_T^*),$$

where $c_6$ does not depend on $T$, $w$, $\psi$, $y$ and $z$. The result now follows from Proposition 1. $\square$

LEMMA 3. *Let the assumptions of Theorem 3 hold. Assume, moreover, that $\theta_i \leq \mu_{ij_0}$ for $i \in \mathcal{I}$. Then the nonidling property holds.*

PROOF. It is more convenient here to work with the equivalent set of equations (33)–(35). Fix $T > 0$ and consider $w_i$ strictly increasing, right-continuous on $[0, T)$, with $w_i(0) > 0$. For $i \neq i_0$, by (34), $\psi_{ij_0} = x_i - y_i$, hence by (33), $x_i = w_i + \int_0^t [-\mu_{ij_0} x_i(s) + (\mu_{ij_0} - \theta_i) y_i(s)] \, ds$. Since the second term in the integrand is positive, by a standard comparison result $x_i \geq \xi_i$, where $\xi_i$ solves $\xi_i = w_i - \mu_{ij_0} \int_0^t \xi_i(s) \, ds$. By Lemma 1, and since $w_i(0) \geq 0$ and $w_i$ is nondecreasing, we have $\xi_i \geq 0$. As a result,

$$x_i(t) \geq 0, \qquad i \neq i_0, \ t < T. \tag{75}$$

Recall that $e \cdot z = (e \cdot x)^-$ and let $\tau = \inf\{t : e \cdot z(t) > 0\} = \inf\{t : e \cdot x(t) < 0\}$. Arguing by contradiction, assume that $\tau < T$. By right-continuity and (33),

(76) $$e \cdot x(\tau) \leq 0.$$

By (33) and the assumptions on $w$, $\tau > 0$. On $[0, \tau)$, $\psi_{i_0 j} = z_j = 0$, $j \neq j_0$. Applying the argument in the previous paragraph to $x_{i_0}$,

$$x_{i_0}(t) = w_{i_0}(t) + \int_0^t [-\mu_{i_0 j_0} x_{i_0}(s) + (\mu_{i_0 j_0} - \theta_{i_0}) y_{i_0}(s)] \, ds, \qquad t < \tau.$$

Hence $x_{i_0} \geq \xi$ on $[0, \tau)$, where now $\xi = w_{i_0} - \mu_{i_0 j_0} \int_0^{\cdot} \xi$, and using Lemma 1,

$$\xi(t) \geq w_{i_0}(t) \exp(-\mu_{i_0 j_0} t) \geq w_{i_0}(0) \exp(-\mu_{i_0 j_0} T) > 0, \qquad t < \tau.$$

By (33) and since $w_{i_0}$ is nondecreasing on $[0, T)$, it follows that, in fact, $x_{i_0}(\tau) > 0$. Along with (75), this stands in contradiction to (76). Therefore, $e \cdot z$ vanishes on $[0, T)$ and by (68), so do $z_i$, $i \in \mathcal{I}$.    □

All results of this section remain valid for trees of diameter 2, as follows on applying them for trees of diameter 3 and choosing $\mu_{ij} = 0$ for appropriate $(i, j)$ and $\theta_i = 0$ for appropriate $i$ (in this section we have not used the assumption that $\mu_{ij} > 0$, but only $\mu_{ij} \geq \theta_i$).

Theorem 3 and Lemma 3 imply the following.

COROLLARY 2. *Under the assumptions of Theorem* 1(ii), *if* (66)–(69) *are satisfied on* $[0, \infty)$ *and* $w$ *is right-continuous, then* (65) *holds, with* $m_4$ *that is independent of* $w$, $\psi$, $y$, $z$ *and* $t$. *Consequently, for any* $m \geq 1$, *any initial condition* $x \in \mathbb{R}^I$ *and any admissible system* $\pi \in \Pi$,

$$E_x^\pi \|X(t)\|^m \leq m_6 (1 + \|x\|)^{m_6} (1 + t)^{m_6}, \qquad t \geq 0,$$

*where* $m_6$ *does not depend on* $\pi$, $x$ *and* $t$.

**6. Case where cost is bounded below and summary of estimates.** This section treats part (iii) of Theorem 1. First it is shown that there exists an admissible system $\pi_0$ under which the state process satisfies a polynomial growth condition in time. Then it is shown that we can consider a subset of $\Pi$ of admissible systems that, in a sense, switch to $\pi_0$ after some time, without losing optimality. The estimates for $\pi_0$ then remain valid in the switched systems.

PROPOSITION 3. *Let Assumption* 1 *hold and consider the system* (33)–(35). *Fix* $(i_0, j_0) \in \mathcal{E}$ *such that* $\mu_{i_0 j_0} \geq \theta_{i_0}$. *Let* $y = (e \cdot x)^+ u$ *and* $z = (e \cdot x)^- v$, *where* $u(t) = e_{i_0}$ *and* $v(t) = e_{j_0}$ *for all* $t$. *Then the estimate* $\|x(t)\| \leq m_7 (1 + t)^{m_7} \|w\|_t^*$ *holds, where* $m_7$ *does not depend on* $w$ *and* $t$.

PROOF. Let $i_0$ be the root. Then $j_0 \in B(i_0)$. It suffices to show that for $k = 0, 1, \ldots$, if $i \in l_{2k}$ and $j \in l_{2k+1}$, then

$$|\psi_{a(j)j}(t)|, |\psi_{ia(i)}(t)|, |x_i(t)| \le c_k(1+t)^{c_k}\|w\|_t^*.$$

This claim is proved by backward induction on $k$.

*Induction base.* The variable $k$ is the smallest number $n$ such that $l_{2n+2} = \varnothing$. If $l_{2k+1}$ is not empty, let $j \in l_{2k+1}$. Then by (35), $\psi_{a(j)j} = -z_j = 0$. Let $i \in l_{2k}$. For $j \in B(i)$, $\psi_{ij} = 0$. Let $j = a(i)$. Then by (33), $x_i = w_i - \mu_{ij}\Im\psi_{ij}$, and by (34) and the assumptions, $\psi_{ij} = x_i - y_i = x_i$. Thus $x_i$ solves $x_i = w_i - \mu_{ij}\Im x_i$ and by Lemma 1, $|x_i|, |\psi_{ij}| \le c_1'|w_i|^*$ for some constant $c_1'$.

*Induction step.* Assuming the claim holds for $k+1$ we show that it holds for $k$. Let $j \in l_{2k+1}$, $j \ne j_0$. Then by (35) and the assumption, $\psi_{a(j)j} = -\sum_{i \in B(j)} \psi_{ij}$. Thus by the induction assumption $|\psi_{a(j)j}| \le c_2'(1+t)^{c_k}\|w\|^*$. Let $i \in l_{2k}$, $i \ne i_0$. Then by (33) and (34), the bound on $\psi_{ij}$, $j \in B(i)$, just obtained, and the induction assumption,

$$x_i(t) = g(t) - \mu_{ia(i)}\Im x_i(t),$$

where $|g(t)| \le c_3'(1+t)^{c_3'}\|w\|_t^*$ with some constant $c_3'$. Hence by Lemma 1, $|x_i(t)| \le c_4'(1+t)^{c_4'}\|w\|_t^*$. By (34) and the induction assumption, a similar bound then holds for $|\psi_{ia(i)}|$.

Consider now $i_0$ and $j_0$. Write $y_0$, $z_0$ and $\psi_0$ for $y_{i_0}$, $z_{j_0}$ and $\psi_{i_0j_0}$, respectively, with a similar convention for $x_0$, $\mu_0$ and $\theta_0$. By (33)–(35) and the induction assumption, we have

$$x_0 = g_1 - \mu_0\Im\psi_0 - \theta_0\Im y_0,$$

$$\psi_0 + g_2 = x_0 - y_0,$$

$$\psi_0 + g_3 = -z_0,$$

where $|g_i(t)| \le c_4'(1+t)^{c_4'}\|w\|_t^*$ and $c_4'$ is independent of $t$ and $w$. Then $x_0 = g_1 + \Im[-\mu_0 g_2 - \mu_0 x_0 + (\mu_4 - \theta_0)y_0]$. Since by assumption the last term in the integral is nonnegative, $x_0 \ge \xi$, where $\xi = g_1 - \mu_0\Im(g_2 + \xi)$. By Lemma 1, $|\xi(t)| \le g_4 = c_5'(1+t)^{c_5'}\|w\|_t^*$ for all $t$. This establishes an appropriate lower bound on $x_0$. Moreover, if $y_0(t) = 0$, then $\psi_0(t) \ge x_0(t) - g_2(t) \ge -g_4(t) - g_2(t)$, and if $y_0(t) > 0$, then $y_0(t) \le e \cdot y(t) = e \cdot x(t)$ and $\psi_0(t) \ge x_0(t) - e \cdot x(t) - g_2(t) = -\sum_{i \ne i_0} x_i(t) - g_2(t)$. Thus, by the induction assumption, $\psi_0(t) \ge -c_6'(1+t)^{c_6'}\|w\|_t^*$ holds for all $t$. Since $\psi_0 \le -g_3$, a similar upper bound on $\psi_0$ holds. Finally, $x_0 \le |g_1(t)| + \mu_0 t|\psi_0|_t^*$ and, therefore, a similar bound holds for $x_0$ as well. This completes the proof by induction. □

Throughout the rest of this section let $i_0$ and $j_0$ be as in Proposition 3. The following result shows that it suffices to consider only a certain subset of the set of admissible systems that have the following property: There is a stopping time $\vartheta$ on the filtration associated with $\pi$, such that on $\{\vartheta < \infty\}$, $u(t) = e_{i_0}$ and $v(t) = e_{j_0}$ for all $t \geq \vartheta$.

Recall that on a given admissible system $\pi$, for every initial data $x$ there is a corresponding controlled process (by Proposition A.1).

PROPOSITION 4. *Let the assumptions of Theorem 1(iii) hold. Let $x \in \mathbb{R}^I$ be given. Then there is a set of admissible systems $\widetilde{\Pi} \subset \Pi$ such that the following statements are valid.*

(i) *For every $m \geq 1$,*

$$E_x^\pi \|X(t)\|^m \leq m_8 (1+t)^{m_8}, \qquad \pi \in \widetilde{\Pi},$$

*where $m_8$ depends on $x$, $\pi$ and $m$, but not on $t$.*

(ii) *We have $V(x) = \widetilde{V}(x)$, where $\widetilde{V}(x) = \inf_{\pi \in \widetilde{\Pi}} C(x, \pi)$.*

PROOF. For any $\pi \in \Pi$ and any stopping time $\vartheta$ on $\pi$, let $\pi^\vartheta$ be the admissible system obtained from $\pi$ by setting $u(t) = e_{i_0}$ and $v(t) = e_{j_0}$ for $t \geq \vartheta$. Given $x$ and $\pi \in \Pi$ such that $C(x, \pi) < \infty$ and for $\varepsilon > 0$, let

$$(77) \qquad \sigma_\varepsilon = \inf\left\{ t : E_x^\pi \int_t^\infty e^{-\gamma s} L(X(s), U(s)) \, ds \leq \varepsilon \right\}$$

and define the stopping time

$$(78) \qquad \vartheta_\varepsilon = \inf\{ t \geq \sigma_\varepsilon : \|X(t)\| \geq \varepsilon^{-1} \}.$$

Note that $\sigma_\varepsilon$ and $\vartheta_\varepsilon$ depend on $x$ and $\pi$. Note also that

$$(79) \qquad \vartheta_\varepsilon \geq \sigma_\varepsilon \to \infty \qquad \text{as } \varepsilon \to 0.$$

For short, write $\pi^\varepsilon$ for $\pi^{\vartheta_\varepsilon}$. Define

$$\widetilde{\Pi} = \{\pi^\varepsilon : \pi \in \Pi, \varepsilon \in (0,1)\}.$$

For part (i), let $\pi \in \Pi$ and $\varepsilon$ be given. Then on $\pi^\varepsilon$ and on the event $\{\vartheta_\varepsilon = \infty\}$, $\sup_{s \in [0,\infty)} \|X(s)\| \leq \varepsilon^{-1} \vee \|X\|_{\sigma_\varepsilon}^*$. Moreover, by Proposition 3, on the event $\{\vartheta_\varepsilon < \infty\}$, for $t \geq \sigma_\varepsilon$,

$$\|X(t)\| \leq m_7 (1 + t - \sigma_\varepsilon)^{m_7} \left( \|X\|_{\sigma_\varepsilon}^* + \sup_{u \in [\sigma_\varepsilon, t]} \|\widetilde{W}(u) - \widetilde{W}(\sigma_\varepsilon)\| \right).$$

Recall that $\widetilde{W}$ is a Brownian motion with drift. Hence $E_x^\pi (\sup_{u \in [\sigma_\varepsilon, t]} \|\widetilde{W}(u) - \widetilde{W}(\vartheta_\varepsilon)\|) \leq c_1 (1 + t - \sigma_\varepsilon)^{c_1}$ for some $c_1$ independent of $\varepsilon$ and $t$, and the same holds for expectation under $E_x^{\pi^\varepsilon}$. Hence it suffices to prove that $E_x^{\pi^\varepsilon} [(\|X\|_{\sigma_\varepsilon}^*)^m] = E_x^\pi [(\|X\|_{\sigma_\varepsilon}^*)^m] < \infty$. This follows from an easy application of Gronwall's lemma, using the fact that $x \mapsto b(x, U)$ is Lipschitz uniformly in $(x, U)$. Part (i) follows.

By Assumption [2](iv) and Proposition [3],

$$E_x^{\pi^\varepsilon}\left[\mathbb{1}_{\vartheta_\varepsilon<\infty}\int_{\vartheta_\varepsilon}^\infty e^{-\gamma t}L(X(t),U(t))\,dt\right]$$

$$\leq c_2 E_x^{\pi^\varepsilon}\left[\mathbb{1}_{\vartheta_\varepsilon<\infty}\int_{\vartheta_\varepsilon}^\infty e^{-\gamma t}(1+t-\vartheta_\varepsilon)^{c_2}\right.$$

$$\left.\times(1+\|X(\vartheta_\varepsilon)\|+\|\widetilde{W}_t-\widetilde{W}_{\vartheta_\varepsilon}\|)^{m_L}\,dt\right]$$

(80)

$$\leq c_3\int_0^\infty e^{-\gamma t}(1+t)^{c_2}\,dt\,E_x^{\pi^\varepsilon}[\mathbb{1}_{\vartheta_\varepsilon<\infty}e^{-\gamma\vartheta_\varepsilon}(1+\|X(\vartheta_\varepsilon)\|)^{m_L}]$$

$$+c_3 E_x^{\pi^\varepsilon}\left[\mathbb{1}_{\vartheta_\varepsilon<\infty}\int_{\vartheta_\varepsilon}^\infty e^{-\gamma t}(1+t-\vartheta_\varepsilon)^{c_2}\|\widetilde{W}_t-\widetilde{W}_{\vartheta_\varepsilon}\|^{m_L}\,dt\right]$$

$$:=c_4 E_x^\pi[\mathbb{1}_{\vartheta_\varepsilon<\infty}(1+\|X(\vartheta_\varepsilon)\|)^{m_L}e^{-\gamma\vartheta_\varepsilon}]+\alpha_2(\varepsilon)$$

$$:=\alpha_1(\varepsilon)+\alpha_2(\varepsilon),$$

where $c_2$, $c_3$ and $c_4$ do not depend on $\varepsilon$. Conditioning on $F_{\vartheta_\varepsilon}$, using strong Markov and stationary increments properties of $\widetilde{W}$, and using the fact $\vartheta_\varepsilon\geq\sigma_\varepsilon$, we have

$$\alpha_2(\varepsilon)\leq c_3 E_x^{\pi^\varepsilon}\exp(-\gamma\sigma_\varepsilon)E_x^{\pi^\varepsilon}\int_0^\infty e^{-\gamma s}(1+s)^{c_2}\|\widetilde{W}_s\|^{m_L}\,ds$$

$$\leq c_5\exp(-\gamma\sigma_\varepsilon),$$

where $c_5<\infty$ by properties of Brownian motion. Therefore, by [(79)], $\alpha_2$ converges to zero as $\varepsilon\to0$.

Next we show that $\alpha_1(\varepsilon)\to0$ as $\varepsilon\to0$. Below, we sometimes write $\vartheta$ for $\vartheta_\varepsilon$. By definition of $b$ [see [(27)]], $\|b(x,U)\|\leq c_6(1+\|x\|)$ and $\|b(x,U)-b(y,U)\|\leq c_6\|x-y\|$, where $c_6$ does not depend on $x$, $y$ and $U$. Thus by [(28)], for any $t\geq\vartheta$, we have on the event $\{\vartheta<\infty\}$,

$$\|X(t)-X(\vartheta)\|$$

$$\leq c_7\|W(t)-W(\vartheta)\|$$

$$+c_7\int_\vartheta^t[b(X(\vartheta),U(s))+(b(X(s),U(s))-b(X(\vartheta),U(s)))]\,ds$$

$$\leq c_8\|W(t)-W(\vartheta)\|+c_8(t-\vartheta)(1+\|X(\vartheta)\|)+c_8\int_\vartheta^t\|X(s)-X(\vartheta)\|\,ds,$$

where $c_7$ and $c_8$ not depend on $t$. Hence by Gronwall's lemma,

(81)
$$\|X(t)-X(\vartheta)\|$$

$$\leq c_8\left(\sup_{s\in[\vartheta,t]}\|W(s)-W(\vartheta)\|+(t-\vartheta)(1+\|X(\vartheta)\|)\right)\exp(c_8(t-\vartheta)).$$

Let $\tau = \inf\{t > \vartheta : \|X(t)\| \le \|X(\vartheta)\|/2\}$ (and $\tau = \infty$ on $\{\vartheta = \infty\}$). By (78), $\|X(\vartheta)\| \ge \varepsilon^{-1}$. Hence by (81), assuming $\varepsilon$ is small enough, we have on $\{\vartheta < \infty\}$

$$P(\tau - \vartheta < \varepsilon^{1/2} | F_\vartheta)$$

$$\le P\Big( \sup_{t \in [\vartheta, \vartheta + \varepsilon^{1/2}]} \|X(t) - X(\vartheta)\| \ge (1 + \|X(\vartheta)\|)/3 \Big| F_\vartheta \Big)$$

$$\le P\Big( \sup_{t \in [\vartheta, \vartheta + \varepsilon^{1/2}]} c_8 \|W(t) - W(\vartheta)\|$$

(82)

$$\ge (1 + \|X(\vartheta)\|)(3^{-1} - c_8 \exp(c_8 \varepsilon^{1/2}) \varepsilon^{1/2}) \Big| F_\vartheta \Big)$$

$$\le P\Big( \sup_{s \in [\vartheta, \vartheta + \varepsilon^{1/2}]} c_8 \|W(s) - W(\vartheta)\| \ge (6\varepsilon)^{-1} \Big| F_\vartheta \Big)$$

$$\le c_9 \exp(-\varepsilon^{-1}).$$

Denote $\beta_\varepsilon = \tau_\varepsilon \wedge [\vartheta_\varepsilon + \varepsilon^{1/2}]$. Then by (77), (78) and (82), and using the lower bound on $L(x, U)$, for all $\varepsilon$ small,

$$\varepsilon \ge E_x^\pi \Big[ \mathbb{1}_{\vartheta_\varepsilon < \infty} \int_{\vartheta_\varepsilon}^\infty e^{-\gamma s} L(X(s), U(s)) \, ds \Big]$$

$$\ge c_{10} E_x^\pi \Big[ \mathbb{1}_{\vartheta_\varepsilon < \infty} \int_{\vartheta_\varepsilon}^{\beta_\varepsilon} e^{-\gamma s} (1 + \|X(s)\|)^{m_L} \, ds \Big]$$

$$\ge c_{10} \exp(-\gamma \varepsilon^{1/2}) E_x^\pi [\mathbb{1}_{\vartheta_\varepsilon < \infty} (1 + (1/2) \|X(\vartheta_\varepsilon)\|)^{m_L} \exp(-\gamma \vartheta_\varepsilon)(\beta_\varepsilon - \vartheta_\varepsilon)]$$

$$\ge c_{11} E_x^\pi \{ \mathbb{1}_{\vartheta_\varepsilon < \infty} (1 + (1/2) \|X(\vartheta_\varepsilon)\|)^{m_L}$$

$$\times \exp(-\gamma \vartheta_\varepsilon) \varepsilon^{1/2} P_x^\pi [\tau_\varepsilon - \vartheta_\varepsilon \ge \varepsilon^{1/2} | F_{\vartheta_\varepsilon}] \}$$

$$\ge c_{12} \varepsilon^{1/2} E_x^\pi \{ \mathbb{1}_{\vartheta_\varepsilon < \infty} (1 + \|X(\vartheta_\varepsilon)\|)^{m_L} e^{-\gamma \vartheta_\varepsilon} \},$$

where the constants $c_{10}, c_{11}, c_{12} > 0$ do not depend on $\varepsilon$. As a result, $\alpha_1(\varepsilon) \to 0$ as $\varepsilon \to 0$. Thus by (80), for every $\pi \in \Pi$ and $\varepsilon$ small enough, $\pi^\varepsilon \in \tilde{\Pi}$ satisfies

(83)                    $$C(x, \pi^\varepsilon) \le C(x, \pi) + \alpha(\varepsilon),$$

and $\alpha(\varepsilon) = \alpha_1(\varepsilon) + \alpha_2(\varepsilon) \to 0$ as $\varepsilon \to 0$. Hence $\tilde{V}(x) \le V(x)$ and part (ii) of the result follows.  $\square$

The following proposition summarizes our estimates in cases (i)–(iii) of Theorem 1.

PROPOSITION 5. *In cases* (i)–(iii) *of Theorem* 1, *we have:*

(i) *For any $x$, any admissible system $\pi \in \Pi$ [$\pi \in \widetilde{\Pi}$ in case* (iii)] *and $m \geq 1$,*

$$E_x^\pi \|X(t)\|^m \leq m_9(1+t)^{m_9}, \qquad t \geq 0,$$

*where the constant $m_9$ does not depend on $t$ (but may depend on $x$, $\pi$ and $m$).*

(ii) *There is a constant $m_{10}$, not depending on $x$, such that $V(x) \leq m_{10}(1 + \|x\|)^{m_{10}}$.*

PROOF. Item (i) follows from Corollaries 1 and 2 and Proposition 4, respectively. Item (ii) follows from Proposition 3 [in cases (i) and (ii) alternatively from Corollaries 1 and 2]. □

## APPENDIX

PROPOSITION A.1. *Let initial data $x \in \mathbb{R}^I$ and an admissible system $\pi \in \Pi$ be given. Then there exists a controlled process $X$ associated with $x$ and $\pi$. Moreover, if $X$ and $\overline{X}$ are controlled processes associated with $x$ and $\pi$, then $X(t) = \overline{X}(t)$, $t \geq 0$, $P$-a.s.*

PROOF. Note that $(x, U) \mapsto b(x, U)$ is continuous and $x \mapsto b(x, U)$ is Lipschitz uniformly in $U$. Consider $b_m$, a function that agrees with $b$ on the ball $B(0, m)$, and is uniformly Lipschitz and bounded. Then strong existence and uniqueness for

$$X_m(t) = x + rW(t) + \int_0^t b_m(X_m(s), U(s)) \, ds, \qquad 0 \leq t < \infty,$$

holds by Theorem I.1.1 of [5]. Since $\|X_m(t)\| \leq \|x\| + c\|W(t)\| + c\int_0^t \|X_m(s)\| \, ds$, we have $\|X_m(t)\| \leq (\|x\| + c\|W\|_t^*)(1 + e^{ct})$ by Gronwall's lemma. Thus letting $\tau_m = \inf\{t : \|X_m(t)\| \geq m\}$, we have $\tau_m \to \infty$ a.s. Therefore, $X(t) = \lim_m X_m(t)$ for all $t$ defines a process that solves the equation (a strong solution). If $X$ and $\overline{X}$ are both strong solutions, then for every $m$, they both agree with $X_m$ on $[0, \tau_m]$. Therefore, they agree on $[0, \infty)$ a.s. □

PROPOSITION A.2. *Let Assumption 1 hold. Then given $\alpha_i, \beta_j \in \mathbb{R}$, $i \in \mathcal{I}$, $j \in \mathcal{J}$, satisfying $\sum \alpha_i = \sum \beta_j$, there exists a unique solution $\psi_{ij}$ to the set of equations*

$$(84) \qquad \sum_{j \in \mathcal{J}} \psi_{ij} = \alpha_i, \qquad i \in \mathcal{I},$$

$$(85) \qquad \sum_{i \in \mathcal{I}} \psi_{ij} = \beta_j, \qquad j \in \mathcal{J},$$

*where $\psi_{i,j} = 0$ for $i \not\sim j$.*

PROOF.    We use notation from Section 4. Let $i_0 \in \mathcal{I}$ be the root. We show the following claim by backward induction on $k$: For $k \in [1, K]$ even (resp. odd), if $i \in l_k$ (resp. $j \in l_k$), then $\psi_{ia(i)}$ (resp. $\psi_{ja(j)}$) is uniquely determined by (84) and (85).

*Induction base.* The variable $k$ is the largest $n$ such that $l_n$ is nonempty. If $k$ is even, let $i \in l_k$. Then $B(i)$ is empty and (84) implies $\psi_{ia(i)} = \alpha_i$. The case $k$ odd is similar.

*Induction step.* Assume the claim holds for $k$. Consider the case where $k$ is odd (the case $k$ even is treated similarly). For $i \in l_{k-1}$, (84) shows $\psi_{ia(i)} = \alpha_i - \sum_{j \in B(i)} \psi_{ij}$, and since by the induction assumption, $\psi_{ij}$ are uniquely determined for $j \in B(i)$, so is $\psi_{ia(i)}$. This completes the proof by induction and the result follows.    $\square$

PROOF OF THEOREM 1.    Based on Proposition 5, the proof of Theorem 1 is similar to the proof of Theorem 2 of [4]. Since this is the main result of this paper, we repeat it here with modifications, mainly to accommodate case (iii).

We first consider (29) on a smooth open bounded connected domain $\Gamma$, satisfying an exterior sphere condition, with boundary conditions

(86)                    $$f(x) = V(x), \qquad x \in \partial\Gamma.$$

The key is a result from [9] regarding existence of classical solutions in bounded domains with merely continuous boundary conditions. To use this result, we verify the following two conditions.

  (i)  We have $|H(x, p)| \leq c(1 + \|p\|)$ for $x \in \Gamma$, where $c$ does not depend on $x$ or $p$.

  (ii)  We have $H(x, p) \in C^\varepsilon(\overline{\Gamma} \times \mathbb{R}^I)$ for some $\varepsilon \in (0, 1)$.

Item (i) is immediate from the local boundedness of $b(x, U)$ and $L(x, U)$. Next we show that item (ii) holds. For $\delta > 0$, let $V$ be such that $H(y, q) \geq b(y, V) \cdot q + L(y, V) - \delta$. Write

$$H(x, p) - H(y, q) \leq b(x, V) \cdot p + L(x, V) - b(y, V) \cdot q - L(y, V) + \delta.$$

Using the Hölder property of $L$ in $x$ uniformly for $(x, V) \in \overline{\Gamma} \times \mathbb{U}$ and the Lipschitz property of $b$ in $x$ uniformly in $(x, V)$,

$$H(x, p) - H(y, q) \leq c\|p - q\| + c\|p\| \|x - y\| + c\|x - y\|^\rho + \delta.$$

Since $\delta > 0$ is arbitrary, it can be dropped. This shows that $H$ is Hölder continuous with exponent $\rho$, uniformly over compact subsets of $\overline{\Gamma} \times \mathbb{R}^I$. Hence (ii) holds.

Defining for $(x, z, p) \in \Gamma \times \mathbb{R} \times \mathbb{R}^I$, $A(x, z, p) = (1/2)r^2 p$ and $B(x, z, p) = H(x, p) - \gamma z$, we can write (29) in divergence form as

$$\mathrm{div} A(x, f, Df) + B(x, f, Df) = 0.$$

The hypotheses of Theorem 15.19 of [9] regarding the coefficients $A$ and $B$ hold in view of (i) and (ii). Indeed, $B$ is Hölder continuous of exponent $\rho$, uniformly on compact subsets of $\Gamma \times \mathbb{R} \times \mathbb{R}^I$. Moreover, with $\tau = 0$, $\nu(z) = (1/2)\min_i r_i^2$, $\mu(z) = c(1 + \|z\|)$, $\alpha = 2$, $b_1 = 0$ and $a_1 = 0$, we check that the conditions (15.59), (15.64), (15.66) and (10.23) of [9] are satisfied. Theorem 15.19 of [9] therefore applies [with condition (15.59) instead of (15.60)]. It states that there exists a solution to (29) in $C^{2,\rho}(\Gamma) \cap C(\overline{\Gamma})$, satisfying the continuous boundary condition (86). We denote this solution by $f$.

Let $x \in \Gamma$. Let $\pi$ be any admissible system in $\Pi$ and let $X$ be the controlled process associated with $x$ and $\pi$. Let $\tau$ denote the first time $X$ hits $\partial\Gamma$. Using Itô's formula for the $C^{1,2}(\mathbb{R}_+ \times \Gamma)$ function $e^{-\gamma t}f(x)$, in conjunction with the inequality

$$\mathcal{L}f(y) + b(y,U) \cdot Df(y) + L(y,U) - \gamma f(y) \geq 0, \qquad y \in \Gamma, U \in \mathbb{U},$$

satisfied by $f$, we obtain

(87)
$$f(x) \leq \int_0^{t \wedge \tau} e^{-\gamma s} L(X_s, U_s)\, ds + e^{-\gamma(t \wedge \tau)} f(X_{t \wedge \tau})$$
$$- \int_0^{t \wedge \tau} e^{-\gamma s} Df(X_s) \cdot r\, dW_s.$$

Taking expectation and then sending $t \to \infty$, using monotone convergence for the first term and bounded convergence for the second term, we have with $g(t,x) = e^{-\gamma t}V(x)$,

$$f(x) \leq E_x^\pi\left[\int_0^\tau e^{-\gamma s} L(X_s, u_s)\, ds + e^{-\gamma\tau} V(X_\tau)\right] = C_{\Gamma,g}(x,\pi).$$

Taking the infimum over $\pi \in \Pi$, we have $f(x) \leq V_{\Gamma,g}(x) = V(x)$, $x \in \Gamma$, by Proposition 5(iv).

To obtain the equality $f = V$ on $\Gamma$, we next show there exist optimal Markov control policies for the control problem on $\Gamma$. Let

(88)    $$\varphi(x,U) = b(x,U) \cdot Df(x) + L(x,U), \qquad x \in \Gamma, U \in \mathbb{U}.$$

Note that $\varphi$ is continuous on $\Gamma \times \mathbb{S}^k$. For each $x$, consider the set $\mathbf{U}_x \neq \varnothing$ of $U \in \mathbb{U}$ for which

$$\varphi(x,U) = \inf_{V \in \mathbb{U}} \varphi(x,V).$$

We show that there exists a measurable selection of $\mathbf{U}_x$, namely there is a measurable function $h$ from $(\Gamma, \mathcal{B}(\Gamma))$ to $(\mathbb{U}, \mathcal{B}(\mathbb{U}))$ with $h(x) \in \mathbf{U}_x$, $x \in \Gamma$.

Let $x_n \in \Gamma$ and assume $\lim_n x_n = x \in \Gamma$. Let $U_n$ be any sequence such that $U_n \in \mathbf{U}_{x_n}$. We claim that any accumulation point of $U_n$ is in $\mathbf{U}_x$. For example, if this is not true, then by continuity of $\varphi$, there is a converging subsequence $U_m$, converging to $\widetilde{U}$, and there is a $\widehat{U}$ such that $\delta := \varphi(x, \widehat{U}) -$ $\varphi(x, \widehat{U}) > 0$. Hence for all $m$ large, $\varphi(x_m, U_m) \geq \varphi(x, \widehat{U}) + \delta/2 \geq \varphi(x_m, \widehat{U}) + \delta/4$, contradicting $U_m \in \mathbf{U}_{x_m}$.

As a consequence, the assumptions of Corollary 10.3 in the Appendix of [6] are satisfied and it follows that there exists a measurable selection $h : \Gamma \to \mathbb{U}$ of $(\mathbf{U}_x, x \in \Gamma)$.

We extend $h$ to $\mathbb{R}^I$ in a measurable way so that it takes values in $\mathbb{U}$ (but otherwise arbitrary). Clearly, $x \mapsto b(x, h(x))$ is measurable. Consider the autonomous stochastic differential equation

$$(89) \qquad X(t) = x + rW(t) + \int_0^t \hat{b}(X_s) \, ds,$$

where $\hat{b}(y)$ agrees with $b(y, h(y))$ on $\Gamma$ and is set to zero off $\Gamma$. Then $\hat{b}$ is measurable and bounded on $\mathbb{R}^I$. By Proposition 5.3.6 of [14], there exists a weak solution to this equation. That is, there exists a complete filtered probability space on which $X$ is adapted and $W$ is an $I$-dimensional Brownian motion, such that (89) holds for $t \geq 0$ a.s. On this probability space, consider the process $U_s = h(X_s)$. Since $X$ has continuous sample paths and is adapted, it is progressively measurable (see Proposition 1.13 of [14]) and by measurability of $h$, so is $U$. Denote by $\pi$ the admissible system thus constructed. Then for $s < \tau$, $U_s \in \mathbf{U}_{X_s}$ and $b(X_s, U_s) \cdot Df(X_s) + L(X_s, U_s) = H(X_s, Df(X_s))$. Hence

$$(90) \qquad \mathcal{L}f(X) + b(X_s, U_s) \cdot Df(X_s) + L(X_s, U_s) - \gamma f(X) = 0, \qquad s < \tau.$$

A use of Itô's formula and the convergence theorems just as before now shows that (87) holds with equality, and

$$f(x) = E_x^\pi \left[ \int_0^\tau e^{-\gamma s} L(X_s, U_s) \, ds + e^{-\gamma \tau} V(X_\tau) \right] = C_{\Gamma, g}(x, \pi), \qquad x \in \Gamma,$$

with $g$ as above. This, together with Proposition 5(iv) shows that $f \geq V_{\Gamma, g} = V$ on $\Gamma$. Summarizing, $f = V$ on $\Gamma$. In particular, $V \in C^{2, \rho}(\Gamma)$ and is a classical solution to the HJB equation. Now $\Gamma$ can be taken arbitrarily large, and this shows that $V \in C^{2, \rho}(\mathbb{R}^I)$ and that it satisfies the HJB equation on $\mathbb{R}^I$. In view of Proposition 5(ii), it also satisfies the polynomial growth condition in cases (i)–(iii) of the main result; in case (iv) it is trivially a bounded function. As a result, there exists a classical solution to (29) in $C^{2, \rho}(\mathbb{R}^I)$, again denoted by $f$, satisfying (30) and, moreover, $V = f$.

It remains to show uniqueness in the appropriate class and existence of optimal Markov control policies for the problem on $\mathbb{R}^I$. In cases (i)–(iii) [resp. (iv)], let $\bar{f} \in C_{\text{pol}}^2(\mathbb{R}^I)$ [resp. $C_b^2(\mathbb{R}^I)$] be a solution to (29). Then analogously to (87), we obtain

$$(91) \qquad \bar{f}(x) \leq \int_0^t e^{-\gamma s} L(X_s, U_s) \, ds + e^{-\gamma t} \bar{f}(X_t) - \int_0^t e^{-\gamma s} D\bar{f}(X_s) \cdot r \, dW_s.$$

Consider first cases (i)–(iii). Taking expectation, sending $t \to \infty$, using the polynomial growth of $\bar{f}$ and the moment bounds on $\|X_t\|$ asserted in Proposition 5(i), we have that $\bar{f}(x) \leq C(x, \pi)$, where $\pi \in \Pi$ [$\pi \in \tilde{\Pi}$ in case (iii)] is arbitrary. In case (iv), the same conclusion holds since $\bar{f}$ is bounded. Consequently, $\bar{f} \leq V$ on $\mathbb{R}^d$.

In cases (i) and (ii), the proof of existence of optimal Markov policies as well as the inequality $V \leq \bar{f}$ on $\mathbb{R}^I$ is completely analogous to that on $\Gamma$, where we replace $\Gamma$ with $\mathbb{R}^I$. The weak existence of solutions to (89) follows on noting that $\hat{b}$ satisfies a linear growth condition of the form $\|\hat{b}(y)\| \leq x(1 + \|y\|)$, $y \in \mathbb{R}^I$, and using again Proposition 5.3.6 of [14]. Then as before, (91) is satisfied with equality, and taking expectation and using the polynomial growth condition of $\bar{f}$ and the moment estimates on $\|X\|$ shows that $V = \bar{f}$ on $\mathbb{R}^I$. We conclude that $f$ is the unique solution in $C^2_{\mathrm{pol}}(\mathbb{R}^I)$ that $V = f$, and that there exists a Markov control policy that is optimal for all $x \in \mathbb{R}^I$. In case (iv) an analogous result is obtained [with uniqueness in $C^2_b(\mathbb{R}^I)$] using the boundedness of $\bar{f}$.

Finally, in case (iii) there is no guarantee that the admissible system $\pi$ constructed using (89) is in $\tilde{\Pi}$ and, therefore, the term $e^{-\gamma t} E^\pi_x \bar{f}(X_t)$ in (91) (that is satisfied with equality) may not tend to zero as $t \to \infty$. However, in this case we claim only uniqueness among nonnegative functions $\bar{f}$ and, therefore, using Itô's formula and (90) gives

$$\bar{f}(x) = E^\pi_x \left[ \int_0^t e^{-\gamma s} L(X_s, U_s)\, ds + e^{-\gamma s} \bar{f}(X_t) \right] \geq E^\pi_x \left[ \int_0^t e^{-\gamma s} L(X_s, U_s)\, ds \right]$$

and $\bar{f}(x) \geq C(x, \pi) \geq V(x)$. $\quad\square$

**Acknowledgments.** I thank Avi Mandelbaum and Marty Reiman for very valuable discussions, and an Associate Editor and two referees for their comments that significantly improved the exposition.

DEPARTMENT OF ELECTRICAL ENGINEERING
TECHNION—ISRAEL INSTITUTE OF TECHNOLOGY
HAIFA 32000
ISRAEL
E-MAIL: atar@ee.technion.ac.il